\documentclass[final]{siamart190516}

\usepackage{algorithmic}
\usepackage{amsfonts}
\usepackage{booktabs}
\usepackage{enumitem}
\usepackage{forloop}
\usepackage{graphicx}
\usepackage{mathrsfs}
\usepackage{mathtools}
\usepackage{multirow}
\usepackage{multicol}
\usepackage[caption=false]{subfig}
\usepackage{tabularx}
\usepackage{tikz}
\usepackage{ulem}

\crefname{algorithm}{Algorithm}{Algorithms}
\Crefname{algorithm}{Framework}{Frameworks}
\Crefname{ALC@unique}{Line}{Lines}
\crefname{ALC@unique}{line}{lines}
\newcounter{ct}\forloop{ct}{1}%
{\value{ct}<10}
{\crefname{enum\roman{ct}}{part}{parts}}

\newcommand{\trace}{\mathrm{Tr}}

\newcommand{\st}{\mathrm{s.}~\mathrm{t.}}

\newcommand{\Diag}{\mathrm{Diag}}

\newcommand{\one}{\mathbf{1}}

\newcommand{\lrbrace}[1]{\left\{#1\right\}}

\newcommand{\inner}[1]{\left\langle#1\right\rangle}

\newcommand{\snorm}[1]{\Vert#1\Vert}
\newcommand{\eps}{\varepsilon}
\newcommand{\F}{\mathrm{F}}
\newcommand{\bbar}{\overline}
\newcommand{\T}{\top}
\newcommand{\PP}{\mathscr{P}}
\newcommand{\rr}{\pmb{r}}
\newcommand{\aaa}{\pmb{a}}
\newcommand{\dd}{~{\rm d}}

\newcommand{\calS}{\mathcal{S}}

\newcommand{\abs}[1]{\left|#1\right|}
\newcommand{\R}{\mathbb{R}}
\newcommand{\N}{\mathbb{N}}

\newcolumntype{R}{>{$}r<{$}}
\newcolumntype{C}{>{$}c<{$}}
\newcolumntype{L}{>{$}l<{$}}
\newsiamthm{assumption}{Assumption}
\newsiamthm{claim}{Claim}
\newsiamremark{remark}{Remark}
\DeclareMathOperator*{\argmin}{arg\,min}

\allowdisplaybreaks[4]

\title{A Global Optimization Approach for Multi-Marginal Optimal Transport Problems with Coulomb Cost\thanks{\textbf{Funding}: the work of the second author was supported by the National Natural Science Foundation of China (No. 11971066). The work of the third author was supported in part by the National Natural Science Foundation of China (No. 1212500491, 11971466, 11991021), Key Research Program of Frontier Sciences, Chinese Academy of Sciences (No. ZDBS-LY-7022).}}

\author{
	Yukuan Hu\thanks{State Key Laboratory of Scientific and Engineering Computing, Academy of Mathematics and Systems Science, Chinese Academy of Sciences, and University of Chinese Academy of Sciences, China (email: \email{ykhu@lsec.cc.ac.cn}, \email{liuxin@lsec.cc.ac.cn}).}
	\and Huajie Chen\thanks{School of Mathematical Sciences, Beijing Normal University, China (email: \email{chen.huajie@bnu.edu.cn}).}
	\and Xin Liu\footnotemark[2]~\thanks{Corresponding author.}
}
\date{}

\headers{a global optimization Approach for MMOT}{Yukuan Hu, Huajie Chen and Xin Liu}

\graphicspath{{figures/}} 

\begin{document}
	
	\maketitle
	\begin{abstract}
		In this work, we construct a novel numerical method for solving the multi-marginal optimal transport problems with Coulomb cost. This type of optimal transport problems arises in quantum physics and plays an important role in understanding the strongly correlated quantum systems. With a Monge-like ansatz, the orginal high-dimensional problems are transferred into mathematical programmings with generalized complementarity constraints, and thus the curse of dimensionality is surmounted. However, the latter ones are themselves hard to deal with from both theoretical and practical perspective. Moreover in the presence of nonconvexity, brute-force searching for global solutions becomes prohibitive as the problem size grows large. To this end, we propose a global optimization approach for solving the nonconvex optimization problems, by exploiting an efficient proximal block coordinate descent local solver and an initialization subroutine based on hierarchical grid refinements. We provide numerical simulations on some typical physical systems to show the efficiency of our approach. The results match well with both theoretical predictions and physical intuitions, and give the first visualization of optimal transport maps for some two dimensional systems.
	\end{abstract}
	
	
	
	\section{Introduction}
	\label{sec:introduction}
	
	The aim of this paper is to provide an optimization method for the multi-marginal optimal transport (MMOT) problems \cite{santam15optimal,villani08optimal} arising in many-electron physics \cite{colombo15multimarginal,cotar15infinite,seidl07strictly}. Let $d\in\{1,2,3\}$ be the dimension of system, $\Omega\subseteq\R^d$ be a bounded domain where the electrons are located, $N\in\N$ with $N\geq 2$ be the number of electrons and $\rr_i\in\Omega~(i\in\{1,\ldots,N\})$ be the position of the $i$-th electron. For the many-electron system, the MMOT problem with Coulomb cost reads
	\begin{equation}
		\label{MMOT:Coulomb}
		\begin{array}{cl}
			\displaystyle\min_{\gamma} & \displaystyle\int_{\Omega^{N}}     c\big(\rr_1,\ldots,\rr_N\big)\gamma\big(\rr_1,\ldots,\rr_N\big)\dd\rr_1\cdots\dd\rr_N\\
			\displaystyle\text{subject to}~(\st) & \displaystyle\Pi_i\gamma(\rr) = \frac{1}{N}\rho(\rr),\quad i=1,\ldots,N,\quad\forall~\rr\in\Omega,
		\end{array}
	\end{equation}
	where the cost function $c\big(\rr_1,\ldots,\rr_N\big)$ is determined by the electron-electron Coulomb interaction
	\begin{equation}
		\label{cost:Coulomb}
		c\big(\rr_1,\ldots,\rr_N\big) := \sum_{i<j}
		\frac{1}{|\rr_i-\rr_j|} ,
	\end{equation}
	$\gamma\big(\rr_1,\ldots,\rr_N\big)$ is an $N$-point probability measure on $\Omega^{N}$, with the single-electron density $\rho:\Omega\to\R$ being the $i$-th marginal $\Pi_i\gamma$, i.e. for any $\rr\in\Omega$, $\frac{1}{N}\rho(\rr)$ equals
	\begin{equation}
		\label{marginal}
		\Pi_i\gamma(\rr) := \int_{\Omega^{N-1}} \gamma\big(\rr_1,\ldots,\rr_{i-1},\rr,\rr_{i+1},\ldots,\rr_N\big) \dd\rr_1\cdots\dd\rr_{i-1}\dd\rr_{i+1}\cdots\dd\rr_N .
	\end{equation}
	Note that the Coulomb interaction $1/|\rr_i-\rr_j|$ between the electrons in \eqref{cost:Coulomb} can be approximated or regularized, especially in the simulations of systems with $d<3$ \cite{bednarek03,friesecke21}.
	Nevertheless, the approach constructed in this paper will make no essential difference	as long as the interaction between the electrons is repulsive (i.e. the cost decreases with respect to $|\rr_i-\rr_j|$).
	Therefore in this paper, we will focus ourselves on the Coulomb interaction of the form \eqref{cost:Coulomb}.
	
	The MMOT problem \cref{MMOT:Coulomb} with Coulomb cost \cref{cost:Coulomb} arises as the strictly correlated electrons limit in the density functional theory (DFT). 
	DFT has been most widely used for electronic structure calculations in physics, chemistry, and material sciences (see \cite{becke14perspective} for a review).
	The strictly correlated electrons limit was first introduced in \cite{seidl99strong}, later noticed in \cite{buttazzo12optimal,cotar13density} that the limit problem is an optimal transport problem. The strictly correlated electrons limit provides an alternative route to derive the DFT energy functionals and has been exploited to extend the capability of DFT to treat strongly correlated quantum systems \cite{chen15pair,chen14numerical,grossi17fermionic,malet12strong,mendl14wigner}.
	
	\par Direct discretization of the MMOT problem \cref{MMOT:Coulomb} leads to a linear programming, with the size increasing exponentially fast with respect to $N$ (the number of electrons/marginals). 
	There are several works devoted to numerical methods that try to circumvent the curse of dimensionality. In \cite{benamou16a}, the Sinkhorn scaling algorithm based on iterative Bregman projections was applied to an entropy-regularized discretized MMOT problem of 1D systems. In \cite{buttazzo12optimal,mendl13kantorovich}, the numerical methods based on Kantorovich dual of the MMOT problem were proposed, while there are exponentially many constraints in the dual problem. In \cite{khoo20semidefinite,khoo19convex}, a convex relaxation approach was proposed by imposing certain necessary constraints satisfied by the two-marginal, and the relaxed problem was then solved by semidefinite programming to obtain tight lower bounds for the optimal cost. 
	In \cite{alfonsi2021constrained,alfonsi2021approximation}, the existence of sparse global solutions was established and a constrained overdamped Langevin process was proposed to solve the moment constrained relaxations. 
	In \cite{friesecke21,friesecke18breaking}, the sparsity of optimal solution was rigorously justified and an efficient numerical method was proposed based on column generation and machine learning.
	
	\par The starting point of this work is to approximate the $N$-point measure $\gamma$ by the following ansatz
	\begin{equation}
		\label{ansatz:gammas}
		\gamma(\rr_1,\ldots,\rr_N) = \frac{\rho(\rr_1)}{N}\gamma_2(\rr_1,\rr_2)\cdots\gamma_N(\rr_1,\rr_N) ,
	\end{equation}
	where $\gamma_i:\Omega^2\to\R~(i\in\{2,\ldots,N\})$ satisfies
	\begin{equation}
		\label{conditions:gammas}
		\gamma_i(\rr,\rr')\ge 0,
		\quad
		\int_{\Omega}\gamma_i(\rr,\rr')\dd\rr'=1,
		\quad{\rm and}\quad
		\int_{\Omega}\rho(\rr)\gamma_i(\rr,\rr')\dd\rr=\rho(\rr').
	\end{equation}
	Here we do not have $\gamma_1$ since $\gamma_1(\rr,\rr')=\theta(\rr-\rr')$ by convention, where $\theta$ is the Dirac delta function.
	The condition \cref{conditions:gammas} is derived from the multi-marginal constraints \cref{marginal}.
	From a physical point of view, $\gamma_i(\rr,\rr')$ represents the correlation between the first and the $i$-th electron, which gives the probability of finding the $i$-th electron at $\rr'$ while the first electron is located at $\rr$.
	Under ansatz \cref{ansatz:gammas}, the MMOT problem \cref{MMOT:Coulomb} (with $N>2$) can be rewritten as
	\begin{multline}
		\label{eqn:optimization formulation}
		\displaystyle\min_{\gamma_2,\ldots,\gamma_N}\left\{\sum_{2\le i<j\le N}\int_{\Omega}\int_\Omega\int_\Omega\frac{\rho(\rr)\gamma_i(\rr,\rr')\gamma_j(\rr,\rr'')}{|\rr'-\rr''|} \dd\rr\dd\rr'\dd\rr''\right.\\\left.+\sum_{2\le i\le N}\int_{\Omega}\int_\Omega\frac{\rho(\rr)\gamma_i(\rr,\rr')}{|\rr-\rr'|}\dd\rr\dd\rr' ~:~\gamma_2,\ldots,\gamma_N~{\rm satisfy}~\cref{conditions:gammas}\right\}.
	\end{multline}
	We mention that in the case of $N=2$, the first term in the objective of \eqref{eqn:optimization formulation} vanishes, then the problem is reduced to a linear programming and can be solved by standard algorithms \cite{chen14numerical}. In this work, we focus our attention on the $N\ge3$ settings.
	The formulation \cref{eqn:optimization formulation} amounts to a spectacular dimension reduction, in that the unknowns are $N-1$ transports on $\Omega^2$ instead of the $N$-point measure $\gamma$ on $\Omega^N$. Therefore, the degrees of freedom now scale linearly with respect to $N$ rather than exponentially fast. In particular, the ansatz \cref{ansatz:gammas} includes the Monge state \cite{monge1781,seidl07strictly} by taking $\gamma_i(\rr,\rr')=\theta\big(\rr'-T_i(\rr)\big)$ with $T_i~(i\in\{2,\ldots,N\})$ being the so-called optimal transport map. The Monge formulation gives significant information on the MMOT problem and enjoys physical interpretations; see more discussions in \cref{sec:intro:remarks}.
	
	\par In practical calculations, we need to discretize \cref{eqn:optimization formulation} into some finite dimensional problems. The discretization consists of three steps.
	First, we employ a {\it finite elements} like mesh $\mathcal{T}=\lrbrace{e_k}_{k=1}^K$ to partition the domain $\Omega$ into $K$ non-overlapping elements, i.e. $\bigcup_{k=1}^K e_k=\Omega$ and $e_i\bigcap e_j=\emptyset$ when $i\neq j$. Let $e:=[\abs{e_1},\ldots,\abs{e_K}]^\T\in\R_+^K$ denote the volumes of elements. 
	Second, we approximate the marginal $\rho$ by a vector $\varrho:=[\varrho_1,\ldots,\varrho_K]^\T\in\R_+^K$, where the $k$-th entry $\varrho_k:=\frac{1}{\abs{e_k}}\int_{e_k}\rho(\rr)\dd\rr$ gives the marginal/electron mass on the $k$-th element $e_k$. 
	Finally, the Coulomb interactions and and the transports $\gamma_i~(i=2,\cdots,N)$ can be approximated by the effective interactions and transports between elements, i.e. for any $j,k\in\{1,\ldots,K\}$,
	\begin{equation}
		c_{jk}:=\frac{1}{\abs{e_j}\cdot\abs{e_k}}\int_{e_k}\int_{e_j}\frac{1}{\abs{\rr-\rr'}}\dd\rr\dd\rr'\quad\text{and}\quad
		x_{i,jk}:=\frac{1}{\abs{e_j}\cdot\abs{e_k}}\int_{e_k}\int_{e_j}\gamma_i(\rr,\rr')\dd\rr\dd\rr',
		\label{eqn:integral average}
	\end{equation}
	respectively, leading to $K\times K$ matrices $C:=((1-a_{jk})c_{jk})_{jk}$ and $X_i=(x_{i,jk})_{jk}$ for $i=2,\ldots,N$. Here $a_{jk}$ equals $1$ if $j=k$ and $0$ otherwise. With this discretization, we can approximate \cref{eqn:optimization formulation} using the following optimization problem with unknowns $\{X_i\}_{i=2}^N$:
	\begin{equation}
		\label{eqn:complicate MPCC}
		\begin{array}{cl}
			\displaystyle\min_{X_2,\ldots,X_N} & \displaystyle f(X_2,\ldots,X_N):=\sum_{2\le i\le N}\inner{X_i,\Lambda\Xi C\Xi}+\sum_{2\le i<j\le N}\inner{X_i,\Xi\Lambda X_j\Xi C\Xi}\\
			\st & X_ie=\one,~X_i^\T\Xi\varrho=\varrho,~\trace(X_i)=0,~X_i\ge0,~i=2,\ldots,N,\\
			& \inner{X_i,X_j}=0,~\forall~ i\ne j,
		\end{array}
	\end{equation}
	where $\one$ is the all-one vector in $\R^K$, $\Lambda=\Diag(\varrho),~\Xi=\Diag(e)$ are $K\times K$ diagonal matrices formed by entries in $\varrho$ and $e$, respectively. 
	More detailed derivation of \cref{eqn:complicate MPCC} is given in \cref{appsec:discretization and error analysis}. 
	Note that the diagonal elements in matrix $C$ are removed due to integral divergence in \cref{eqn:integral average}. 
	The extra constraints
	$$\trace(X_i)=0,~i=2,\ldots,N,\quad\text{and}\quad\inner{X_i,X_j}=0,~\forall~ i\ne j$$
	are hence accordingly added. 
	From a physical point of view, this constraint can keep the electrons spatially away from each other in the case of Coulomb repulsion, so that unfavorable particle clustering can be avoided.
	
	\par In the case of $N=3$, \cref{eqn:complicate MPCC} is a mathematical programming with complementarity constraints (MPCC) in view of nonnegative constraints and $\inner{X_2,X_3}=0$. Due to the disjunctive nature of feasible set, a general MPCC violates commonly used constraint qualifications at any feasible point \cite{flegel2005guignard}. As a result, the well-known Karush-Kuhn-Tucker conditions are no longer certificate for feasible points to be local minimizers. When $N>3$, the formulation of the constraints in \cref{eqn:complicate MPCC} is more complicated than that of the complementarity constraints.
	Since $\inner{X_i,X_j},~\forall~ i\ne j$ impose the requirements that, for each $i\in\{2,\ldots,N\}$, the block variable $X_i$ complements all the other blocks, we call \cref{eqn:complicate MPCC} a mathematical programming with generalized complementarity constraints (MPGCC). 
	
	\par In addition to its intrinsic difficulty, we are in quest for global solutions of \cref{eqn:complicate MPCC}. This is a hard matter because both the repulsive energy $f$ and the feasible set are nonconvex in variables $(X_i)_{i=2}^N$. Since the degrees of freedom $(N-1)K^2$ grow quickly as the meshes become finer, the state-of-art global optimization solvers cannot be our last resort.
	
	\subsection{Optimization Background}
	\label{sec:optimization}
	
	\par Although little is known about MPGCC, there exists rich literature on MPCC. To overcome the intrinsic difficulties mentioned above, several MPCC-tailored constraint qualifications have been provided for MPCC. Under these constraint qualifications, points satisfying certain stationary systems are shown to be proper candidates of local minimizers. The related notions and theoretical results are gathered in \cite{scheel2000mathematical,jane2005necessary} and the references within. 
	
	\par With these in place, researchers have proposed various numerical approaches, wherein those based on the original MPCC formulation rank top choices; they employ \textit{the modified nonlinear programming solvers}. For example, the authors in \cite{fletcher2004solving} solved MPCCs using sequential quadratic programming algorithm with filter techniques \cite{fletcher2002nonlinear}; the software introduced in \cite{byrd2006knitro,waltz2006interior} incorporates a suite of nonlinear programming algorithms to tackle MPCCs, including interior-point methods and sequential quadratic programming algorithm, together with globalization techniques such as line search and trust region. 
	
	\par Owing to the troubles when coping with complementarity constraints, methods based on penalty functions gain popularity as well. Among others, we confine our attention to the \textit{$\ell_1$ (complementarity) penalty function}, which favours direct extension to MPGCC \cref{eqn:complicate MPCC} as
	\begin{equation}
		f(X_2,\ldots,X_N)+\beta\sum_{i<j}\inner{X_i,X_j},
		\label{eqn:comp penalty function}
	\end{equation}
	namely, penalizing merely the complementarity violation in $\ell_1$ form. Here $f$ is the repulsive energy defined in \cref{eqn:complicate MPCC}, $\beta>0$ is the penalty parameter.
	Apart from algorithmic benefit, with $N=3$, it can be verified under certain conditions that the global solutions of \cref{eqn:complicate MPCC} coincide with those globally minimizing \cref{eqn:comp penalty function} over $\calS^{N-1}$, where
	\begin{equation}\label{eqn:partial feasible set} \calS:=\{W\in\R^{K\times K}:We=\one,~W^\T\Xi\varrho=\varrho,~\trace(W)=0,~W\ge0\}.\end{equation}
	A direct consequence is that, if the global solutions of \cref{eqn:complicate MPCC} are required, one can in turn minimize \cref{eqn:comp penalty function} over $\calS^{N-1}$ starting with proper initialization. However, we are not aware of any existing method that fully exploits the special structure of \eqref{eqn:comp penalty function}. A customized algorithm is thus needed, particularly in the large-scale context. 
	
	\par In addition, methods based on approximation (smoothing or regularization), augmented Lagrangian functions and full penalization are available as well. We refer interested readers to \cite{Facchinei1999A,hoheisel2013theoretical,hu2004convergence,jia2021augmented,scholtes2001convergence,scholtes1999exact} and the references therein. Compared with methods using modified nonlinear programming solvers or penalty functions, other approaches suffer from an obvious drawback: for a specific MPGCC, the latter ones require solving a sequence of subproblems in the same size to stationarity or even optimality \cite{kanzow2015theprice}. This weakness excludes them from our choices, particularly when the number of grid points $K$ is tremendously large. 
	
	\subsection{Contributions}
	
	Our contributions are three-fold: 
	\begin{enumerate}[label=\textup{(}\arabic*\textup{)}]
		\item \textit{A global optimization approach, equipped with a local solver and a hierarchical initialization subroutine, is constructed for solving \cref{eqn:complicate MPCC}.}
		\vskip 0.05cm
		\par The initialization subroutine (\cref{alg:gr}), derived from hierarchical grid refinements, helps the local solver locate good approximations of global solutions, and hence serves as the core of the proposed global optimization approach (\Cref{frame:ggr}). The proposed approach saves one from brute-force solving large-scale \cref{eqn:complicate MPCC} via plain global optimization methods. Remarkably in \Cref{frame:ggr}, the optimal transport maps can be directly evaluated by the solutions, which is usually difficult in the context of Coulomb cost.
		\vskip 0.1cm
		\item \textit{An inexact proximal block coordinate descent (\textup{\texttt{PBCD}}) algorithm is proposed for locally minimizing \cref{eqn:comp penalty function} over $\calS^{N-1}$.}
		\vskip 0.05cm
		\par The \texttt{PBCD} algorithm (\cref{frame:PBCD scheme}) acts as the local solver in \Cref{frame:ggr} and enjoys global convergence guarantee in the presence of iterate infeasibility (\Cref{thm:rough convergence property}), which is not covered by existing works.
		\vskip 0.1cm
		\item \textit{Simulations of optimal transport maps for some typical 1D and 2D systems.}
		\vskip 0.05cm
		\par We consider systems with the number of electrons up to 7, and discretization with the number of grid points up to $1.6\times10^4$. The results are in line with both theoretical predictions and physical intuitions (\cref{sec:Numerical Experiments}). We also give the first visualization of optimal transport maps for some 2D systems.
	\end{enumerate}

	\subsection{Further Remarks}
	\label{sec:intro:remarks}	
	\mbox{}
	
	\par {\bf Monge formulations.}
	It is unknown whether the MMOT problem \cref{MMOT:Coulomb} with Coulomb cost has a solution of the form \cref{ansatz:gammas}.
	However, the ansatz \cref{ansatz:gammas} includes the Monge solutions, which are most widely studied in physics. The Monge formulation makes the ansatz
	\begin{equation}
		\label{eqn:Monge formulation}
		\gamma(\rr_1,\ldots,\rr_N) = \frac{\rho(\rr_1)}{N}\theta\left(\rr_2-T_2(\rr_1)\right)\cdots \theta\left(\rr_N-T_N(\rr_1)\right) ,
	\end{equation}
	where $\theta$ is the Dirac delta function, the transport map $T_i:\Omega\to\Omega~(i\in\{2,\ldots,N\})$ (we can prescribe $T_1(\rr)=\rr$ for completeness of notations) preserves the single-electron density $\rho$.
	The Monge solution has a simple physical interpretation: the many-electron repulsive energy is minimized at a state such that one electron at position $\rr$ can determine the positions of all other $N-1$ electrons via $\{T_i\}_{i=2}^N$. It is known that for 1D systems, the Monge formulation gives the global solution of the MMOT problems \cite{colombo15multimarginal,cotar13density}. But in the general $d>1$ and $N>2$ cases, it is unknown whether there exists a minimizer of \cref{MMOT:Coulomb} in the form \cref{ansatz:gammas}. Nevertheless, the Monge solution involves a lot of physical information of the many-electron system and can give rise to the Kantorovich potential (which is needed in applications for electronic structure); see \cite{seidl07strictly,seidl17the}. Therefore, the Monge solution is crucial for the MMOT problems in DFT, which is though still difficult to evaluate in the context of Coulomb cost. In our \Cref{frame:ggr}, however, the optimal transport maps $\{T_i\}_{i=2}^N$ can be approximated by the transportation between elements in mesh. More precisely, let $\aaa_j$ be the barycenter of the element $e_j$, then $T_i(\aaa_j)~(i=2,\cdots,N)$ can be approximated by solution $(X_i)_{i=2}^N$ as
	\begin{equation}
		\displaystyle T_i^K(\aaa_j) := \frac{\sum_{1\le k\le K} \aaa_k x_{i,jk}}{\sum_{1\le l\le K}x_{i,jl}},
		\quad j=1,\ldots,K,\quad i=2,\ldots,N .
		\label{eqn:approx transport maps}
	\end{equation}
	
	\vskip 0.2cm
	
	\par {\bf Symmetric constraints.}
	In physics, one is only interested in the measures that are symmetric with respect to $\{\rr_i\}$ (as $\gamma(\rr_1,\ldots,\rr_N)$ represents $N$-point position density of electrons, which is symmetric by the laws of quantum theory). More precisely, one requires that for any permutation $\PP$ on $\{1,\ldots,N\}$, $\gamma\big(\rr_1,\ldots,\rr_N\big) = \gamma\big(\rr_{\PP(1)},\ldots,\rr_{\PP(N)}\big)$. Although we do not have this symmetric restriction in the MMOT problem \cref{MMOT:Coulomb} and the ansatz \cref{eqn:Monge formulation} is in general not symmetric, dropping the restriction does not alter the minimum value. This is because we have a symmetric cost function $c$ in \cref{cost:Coulomb} and equal marginal for any $\Pi_i\gamma$ in \cref{marginal}. Hence each non-symmetric $\gamma$ can give a symmetric one with the same energy value by symmetrization $\frac{1}{N!}\sum_{\PP}\gamma\big(\rr_{\PP(1)},\ldots,\rr_{\PP(N)}\big)$. We do not have to impose the symmetric constraints in the optimization formulation \cref{eqn:complicate MPCC}.
	
	\vskip 0.2cm
	
	\par {\bf Discretization.}
	Most of the existing works discretize the MMOT problems with real space methods \cite{benamou16a,chen14numerical}. Particularly, this paper discretizes \cref{eqn:optimization formulation} into \cref{eqn:complicate MPCC} by representing the marginal $\rho$ with piecewise finite elements and using effective cost coefficients obtained by integrating the continuous cost functions with respect to these elements. To further reduce the computational cost (i.e. use less grid points where the marginal is small), we choose the elements adaptively such that each element carries approximately the same marginal mass.

	\subsection{Outline}
	
	The rest of this paper is organized as follows. We introduce the global optimization approach in \cref{sec:GGR approach}, where the initialization subroutine (\cref{subsec:GR}) and the local solver (\cref{subsec:Algorithm}) are detailed in order. \Cref{sec:Convergence Analysis} is dedicated to the rough statements of the convergence property of \texttt{PBCD}. We corroborate the proposed approach with numerical simulations on several typical systems in \cref{sec:Numerical Experiments}. Finally, conclusions and discussions are drawn in \cref{sec:conclusion}.

	\subsection{Notations}
	
	The image of a linear operator $\mathcal{A}$ is denoted by $\mathrm{Im}(\mathcal{A})$. The notation $\snorm{X}_p$ gives the $p$-norm of matrix $X$, while $\snorm{X}_\F$ yields its Frobenius norm. The components of matrices or vectors are indicated by subscripts, e.g. $x_{ij}$. The inquality $X\ge0$ means $x_{ij}\ge0,~\forall~ i,j$.
	
	\par The notation $\delta_S$ represents the indicator function of set $S$, namely $\delta_S(x)$ equals $0$ if $x\in S$ otherwise $\infty$. For the multi-block objective functions referred in this work (such as \cref{eqn:complicate MPCC}), we occasionally adopt abbreviations in brackets. For example, $f(X_{<i},X_i,X_{>i})$ means $f(X_2,\ldots,X_{i-1},X_i,X_{i+1},\ldots,X_N)$; abbreviations like $X_{<i}$, $X_{(i,j)},X_{>j}$ represent aggregation of blocks with certain subscripts.
	
	\par Regarding algorithm, we use (double) superscripts within bracket for iterates in outer (inner) loop; for instance, $X^{(l)}$ is the iterate in the $l$-th outer iteration, $X^{(l,k)}$ is the iterate in the $k$-th inner iteration of the $(l+1)$-th outer iteration. 
	
	\section{A Global Optimization Approach for Solving (\ref{eqn:complicate MPCC})}
	\label{sec:GGR approach}
	
	In light of ansatz \cref{ansatz:gammas}, the original MMOT problem with Coulomb cost \cref{MMOT:Coulomb} is approximated by MPGCC \cref{eqn:complicate MPCC}. Violating commonly used constraint qualifications, MPGCC \cref{eqn:complicate MPCC} itself is a hard nut to crack from both algorithmic design and theoretical analysis. Rather, we concentrate on the $\ell_1$ penalized MPGCC \cref{eqn:complicate MPCC}, i.e.
	\begin{equation}
		\begin{array}{cl}
			\displaystyle\min_{X_2,\ldots,X_N} & \displaystyle f_\beta(X_2,\ldots,X_N):=f(X_2,\ldots,X_N) + \beta\sum_{i<j}\inner{X_i,X_j}\\
			\st & X_i\in\calS,~i=2,\ldots,N,
		\end{array}
		\label{eqn:penalty problem}
	\end{equation}
	where $f$ is the repulsive energy defined in \cref{eqn:complicate MPCC}, $\calS$, defined in \cref{eqn:partial feasible set}, stands for a section of feasible region. Problem \cref{eqn:penalty problem} is a nonconvex quadratic programming problem, still NP-hard \cite{pardalos1991quadratic}. In the sequel, when we reference \cref{eqn:penalty problem} and its solution in space $\big(\R^{K\times K}\big)^{N-1}$, we simply say \cref{eqn:penalty problem} and its solution with size $K$.
	
	\par For practical purposes, a global solution of \cref{eqn:penalty problem} is always required. Meanwhile, we notice that the degrees of freedom in \cref{eqn:penalty problem}, $(N-1)K^2$, grow fast w.r.t. $K$. This prevents us from brute-force solving \cref{eqn:penalty problem} by state-of-art global optimization methods (e.g. branch-and-bound and cutting plane algorithm) due to exponentially increasing running time. 
	
	\par Motivated by \cite{benamou16a}, we propose a global optimization approach GGR; see \Cref{frame:ggr}.
	\begin{algorithm}
		\floatname{algorithm}{Framework}
		\caption{The GGR approach}
		\label{frame:ggr}
		\begin{algorithmic}[1]
			\REQUIRE{Oracle returning $R,e,\varrho$ in proper dimensions; global solver; local solver; the GR subroutine; initial mesh with $K^{(0)}$ elements $\{e^{(0)}_j\}$.}
			\STATE{Set $l:=0$.}
			\STATE\label{global line}{\textbf{GGR\_Init}: use global solver to solve \cref{eqn:penalty problem} with size $K^{(0)}$ and get $(X_i^{(0)})_{i=2}^N$.}
			\WHILE{certain stopping criteria are not satisfied}
			\STATE\label{mesh refinement}{Refine the last mesh $\{e^{(l)}_j\}$ to $\{e^{(l+1)}_j\}$ with $K^{(l+1)}$ elements.}
			\STATE\label{modify line}{Modify $(X_i^{(l)})_{i=2}^N$ using the GR subroutine to obtain $(X_i^{(l,0)})_{i=2}^N$.}
			\STATE\label{local line}{\textbf{GGR\_LS$(l+1)$}: start local solver from $(X_i^{(l,0)})_{i=2}^N$ to solve \cref{eqn:penalty problem} with size $K^{(l+1)}$ and get $(X_i^{(l+1)})_{i=2}^N$.}
			\STATE{$l:=l+1$.}
			\ENDWHILE
			\RETURN{$(X_i^{(l)})_{i=2}^N\in\big(\R^{K^{(l)}\times K^{(l)}}\big)^{N-1}$.}
		\end{algorithmic}
	\end{algorithm}
	Here, ``G'' and ``GR'' stand for global optimization and the GR initialization subroutine based on hierarchical grid refinement, respectively. \textbf{GGR\_Init} and \textbf{GGR\_LS} are in turn referred to as the initial step invoking a global solver, and the subsequent step invoking a local solver. \Cref{frame:ggr} progresses step by step along with the process of mesh refinements.
	
	\par Justification on the usage of global solver in the initial step (\cref{global line} in \Cref{frame:ggr}) is in order. From the point of applicability, given initial size $K^{(0)}$ of moderate magnitude, globally solving \cref{eqn:penalty problem} is amenable to state-of-art global optimization methods. Considering the necessity, the quality of constructed initial points largely depends on the solutions in the previous step. Hence it is a natural choice for us to invoke a global solver in the initial step. For our choices in implementation, please refer to \cref{subsec:default solver}.
	
	\par Without specification, the mesh refinements (\cref{mesh refinement} in \Cref{frame:ggr}) are done such that the coarse meshes are always embedded into the refined meshes. For more remarks, see \cref{sec:intro:remarks}. Although the refinements are uniform in the numerical simulations of present work (\cref{subsec:1D results,subsec:2D results}), practical implementations focus on the region where marginals vary violently. Nevertheless in the latter circumstances, our GGR approach still works.
	
	\par In what follows, we leave the initialization subroutine part to \cref{subsec:GR} and the local solver part to \cref{subsec:Algorithm}, respectively.
	
	\subsection{Initialization Subroutine based on Grid Refinement}\label{subsec:GR}
	
	Brute-force global optimization of \cref{eqn:penalty problem} becomes impracticable once $K$ grows large. One treatment for this is arming a local solver with good initialization. Roughly speaking, if the energy surface forms a basin around the global solution $(X_i)_{i=2}^N$, the local solver is able to find $(X_i)_{i=2}^N$ provided the initial point lies inside the basin near $(X_i)_{i=2}^N$. This subsection is devoted to the development of the GR subroutine for initialization (\cref{modify line} in \Cref{frame:ggr}). In words, the GR subroutine passes the solution information of previous step on to the current one such that good initialization can be anticipated. Without this process, the located point by local solver is very likely not a global minimizer, resulting in bad solution afterwards.
	
	\par We derive the GR subroutine from some 1D numerical experience: for a particular problem (given oracle of $R,e,\varrho$), the solutions with different sizes share ``similar'' patterns. This phenomenon suggests that we can construct an initial point based on the pattern reflected in the solution with a small size $K$. This point was also observed in \cite{benamou16a}, where the authors supplied a  refinement strategy to meet the accuracy demand with relatively low cost for discretized 1D \cref{MMOT:Coulomb}. Their strategy, however, remains to be explained rigorously and quantitatively. More importantly, they did not discuss the treatment in higher-dimensional context. In the following, we try to understand the ``similarity'' standing at optimal transport and then introduce the GR subroutine. Basically, the proposed subroutine is applicable under any space dimension $d$.
	
	\par
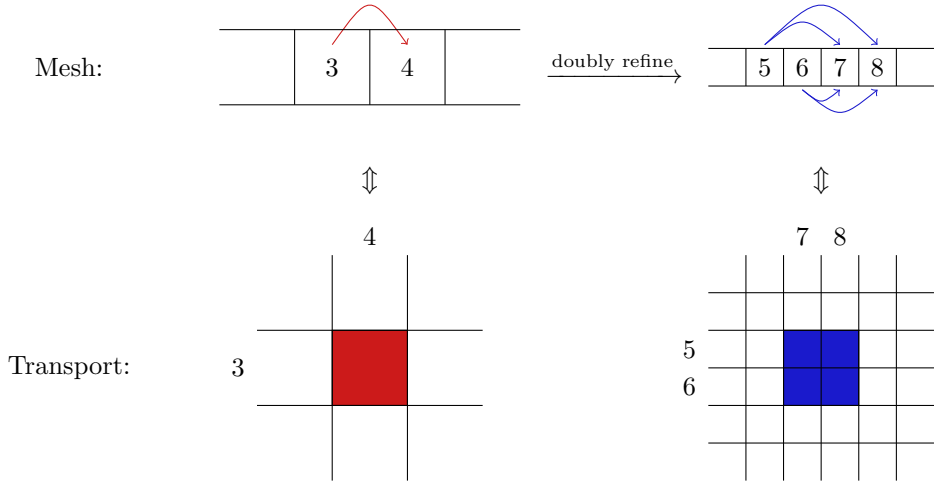
\begin{figure}[htbp]
		\centering
		\begin{tikzpicture}
			\filldraw[fill={rgb:red,8;green,1;blue,1}](1.5,0) -- (2.5,0) -- (2.5,1) -- (1.5,1) -- (1.5,0);
			\filldraw[fill={rgb:red,1;green,1;blue,8}](7.5,0) -- (8.5,0) -- (8.5,1) -- (7.5,1) -- (7.5,0);
			
			\node at (2,2.25) {4};
			\node at (0.25,0.5) {3};
			\node at (7.75,2.25) {7};
			\node at (8.25,2.25) {8};
			\node at (6.25,0.75) {5};
			\node at (6.25,0.25) {6};
			\node at (1.5,4.5) {3};
			\node at (2.5,4.5) {4};
			\node at (7.25,4.5) {5};
			\node at (7.75,4.5) {6};
			\node at (8.25,4.5) {7};
			\node at (8.75,4.5) {8};
			\node at (5.25,4.5) {$\xrightarrow{\text{doubly refine}}$};
			\node at (2,3) {$\Updownarrow$};
			\node at (8,3) {$\Updownarrow$};
			\node at (-2,4.5) {Mesh:};
			\node at (-2,0.5) {Transport:};
			\draw[color={rgb:red,8;green,1;blue,1},->] (1.5,4.8) .. controls (2,5.5) .. (2.5,4.8);
			\draw[color={rgb:red,1;green,1;blue,8},->] (7.25,4.8) .. controls (8,5.5) .. (8.75,4.8);
			\draw[color={rgb:red,1;green,1;blue,8},->] (7.25,4.8) .. controls (7.75,5.2) .. (8.25,4.8);
			\draw[color={rgb:red,1;green,1;blue,8},->] (7.75,4.2) .. controls (8.25,3.8) .. (8.75,4.2);
			\draw[color={rgb:red,1;green,1;blue,8},->] (7.75,4.2) .. controls (8,4.0) .. (8.25,4.2);
			
			\draw[domain=0:4] plot(\x,4);
			\draw[domain=0:4] plot(\x,5);
			\draw[domain=6.5:9.5] plot(\x,4.25);
			\draw[domain=6.5:9.5] plot(\x,4.75);
			
			\foreach \y in {1,2,3}
			\draw[domain=4:5] plot(\y,\x);
			\foreach \y in {7,7.5,8,8.5,9}
			\draw[domain=4.25:4.75] plot(\y,\x);
			
			\foreach \y in {0,1}
			\draw[domain=0.5:3.5] plot(\x,\y);
			\foreach \y in {0,1}
			\draw[domain=-1:2] plot(\y+1.5,\x);
			
			\foreach \y in {-0.5,0,0.5,1,1.5}
			\draw[domain=6.5:9.5] plot(\x,\y);
			\foreach \y in {7,7.5,8,8.5,9}
			\draw[domain=-1.0:2.0] plot(\y,\x);
		\end{tikzpicture}
		\caption{1D case. \textit{The red block means there is mass transported from 3 to 4. Then in a doubly refined mesh, there is mass transported from 5 and 6 to 7 and 8, as is marked out by 4 blue blocks.}}
		\label{fig:gr 1d}
	\end{figure}

 Let us begin with 1D setting. Suppose we already have a finite elements mesh $\{e_j\}$ and a global solution $(X_i)_{i=2}^N$ of \cref{eqn:penalty problem} with $X_i=(x_{i,jk})_{jk}$. Then for any $i\in\{2,\ldots,N\}$, $x_{i,jk}>0$ means that mass of $x_{i,jk}$ is transported from $e_j$ to $e_k$ by transport $X_i$. For the problem with a doubly refined mesh $\{\tilde e_k\}$, the original $e_j,e_k$ correspond to $\tilde e_{2j-1}$ and $\tilde e_{2j}$, $\tilde e_{2k-1}$ and $\tilde e_{2k}$, respectively. Let $j_1=2j-1,~j_2=2j,~k_1=2k-1,~k_2=2k$. A reasonable speculation is that there also exists certain mass transported from $\tilde e_{2j-1},\tilde e_{2j}$ to $\tilde e_{2k-1},\tilde e_{2k}$ by the new $\widetilde X_i=(\tilde x_{i,mn})_{mn}$, i.e. $\tilde x_{i,j_tk_u}>0$, $t,u\in\{1,2\}$, which happens to explain the similarity observed in \cite{benamou16a}. See \Cref{fig:gr 1d} for illustration.
	
	\par
\begin{figure}[htbp]
		\centering
		\scalebox{0.895}{
			\begin{tikzpicture}
				\filldraw[fill={rgb:red,8;green,1;blue,1}](2,1) -- (3,1) -- (3,2) -- (2,2) -- (2,1);
				\filldraw[fill={rgb:red,1;green,1;blue,8}](7.5,1.5) -- (8.5,1.5) -- (8.5,2.5) -- (7.5,2.5) -- (7.5,1.5);
				\filldraw[fill={rgb:red,1;green,1;blue,8}](7.5,-0.5) -- (8.5,-0.5) -- (8.5,0.5) -- (7.5,0.5) -- (7.5,-0.5);
				\filldraw[fill={rgb:red,1;green,1;blue,8}](9.5,1.5) -- (10.5,1.5) -- (10.5,2.5) -- (9.5,2.5) -- (9.5,1.5);
				\filldraw[fill={rgb:red,1;green,1;blue,8}](9.5,-0.5) -- (10.5,-0.5) -- (10.5,0.5) -- (9.5,0.5) -- (9.5,-0.5);
				
				\node at (2.5,3.25) {11};
				\node at (-0.25,1.5) {5};
				\node at (7.75,3.75) {35};
				\node at (8.25,3.75) {36};
				\node at (9.75,3.75) {49};
				\node at (10.25,3.75) {50};
				\node at (6.25,2.25) {9};
				\node at (6.25,1.75) {10};
				\node at (6.25,0.25) {23};
				\node at (6.25,-0.25) {24};
				\node at (1.5,5.5) {(2,4)};
				\node at (2.5,6.5) {(1,5)};
				\node at (8.25,5.75) {\tiny (3,7)};
				\node at (8.25,5.25) {\tiny (4,7)};
				\node at (8.75,5.75) {\tiny (3,8)};
				\node at (8.75,5.25) {\tiny (4,8)};
				\node at (9.25,6.75) {\tiny (1,9)};
				\node at (9.25,6.25) {\tiny (2,9)};
				\node at (9.75,6.75) {\tiny (1,10)};
				\node at (9.75,6.25) {\tiny (2,10)};
				\node at (7.75,1.0) {$\vdots$};
				\node at (8.25,1.0) {$\vdots$};
				\node at (9.75,1.0) {$\vdots$};
				\node at (10.25,1.0) {$\vdots$};
				\node at (9.0,2.25) {$\cdots$};
				\node at (9.0,1.75) {$\cdots$};
				\node at (9.0,0.25) {$\cdots$};
				\node at (9.0,-0.25) {$\cdots$};
				\node at (-2,6) {Mesh:};
				\node at (-2,1.0) {Transport:};
				\node at (2,4.25) {$\Updownarrow$};
				\node at (9,4.25) {$\Updownarrow$};
				\node at (5.5,6) {$\xrightarrow{\text{doubly refine}}$};
				
				\draw[color={rgb:red,8;green,1;blue,1},->] (2.5,6.5) .. controls (2.5,5.5) .. (1.5,5.5);
				\draw[color={rgb:red,1;green,1;blue,8},->] (9.25,6.75) .. controls (8.5,6.5) .. (8.25,5.75);
				\draw[color={rgb:red,1;green,1;blue,8},->] (9.25,6.75) .. controls (8.5,6.5) .. (8.25,5.25);
				\draw[color={rgb:red,1;green,1;blue,8},->] (9.25,6.75) .. controls (8.5,6.5) .. (8.75,5.75);
				\draw[color={rgb:red,1;green,1;blue,8},->] (9.25,6.75) .. controls (8.5,6.5) .. (8.75,5.25);
				
				\foreach \y in {6,7}
				\draw[domain=1:3] plot(\x,\y);
				\foreach \y in {2}
				\draw[domain=5:7] plot(\y,\x);
				
				\foreach \y in {5.5,6,6.5,7}
				\draw[domain=9:11] plot(\x-1,\y);
				\foreach \y in {9.5,10,10.5}
				\draw[domain=5:7] plot(\y-1,\x);
				
				\foreach \y in {-1,0,1}
				\draw[domain=0:4] plot(\x,\y+1.0);
				\foreach \y in {-1,0,1}
				\draw[domain=-2:2] plot(\y+2,\x+1.0);
				
				\foreach \y in {-1.0,-0.5,0,0.5,1.5,2,2.5,3.0}
				\draw[domain=7.5:12.5] plot(\x-1,\y);
				\foreach \y in {-1,-0.5,0,0.5,1.5,2,2.5,3}
				\draw[domain=-1.5:3.5] plot(\y+8.0,\x);
		\end{tikzpicture}}
		\caption{2D case ($7\times7$ rectangular mesh). \textit{The red block means there is mass transported from (1,5) to (2,4). Then in a doubly refined mesh, there is mass transported from (1,9), (1,10), (2,9) and (2,10) to (3,7), (3,8), (4,7) and (4,8), as is marked out by 16 blue blocks.}}
		\label{fig:gr 2d}
	\end{figure}
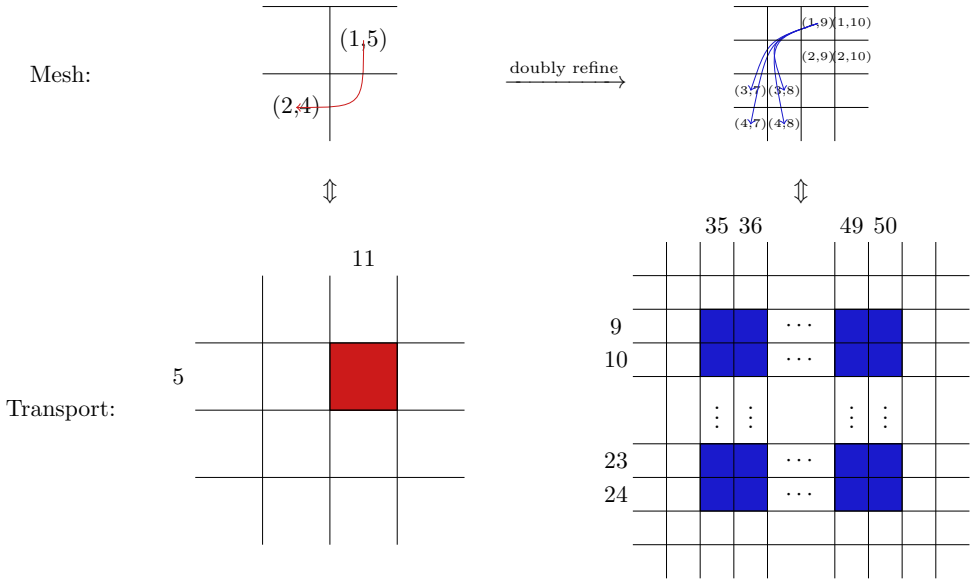

 The above reasoning applies to any $d\in\mathbb{N}$. Suppose a finite elements mesh $\{e_j\}$ and a global solution $(X_i)_{i=2}^N$ are at hand, with $X_i=(x_{i,jk})_{jk}$. Then for any $i\in\{2,\ldots,N\}$, $x_{i,jk}>0$ means that mass of $x_{i,jk}$ is transported from element $e_j$ to $e_k$ by transport $X_i$. After mesh refinement, the original $\{e_j\}$ becomes $\{\tilde e_k\}$; for each $j$, the original element $e_j$ is divided into $s_j$ parts: $e_j=\bigcup_{t=1}^{s_j}\tilde e_{j_t}$ and $\tilde e_{j_{t_1}}\bigcap\tilde e_{j_{t_2}}=\emptyset$ when $j_{t_1}\ne j_{t_2}$. A reasonable speculation is that there also exists certain mass transported from $\tilde e_{j_t}$ to $\tilde e_{k_u}$, where $t\in\{1,\ldots,s_j\},~u\in\{1,\ldots,s_k\}$. Accordingly in $\tilde X_i=(\tilde x_{i,mn})_{mn}$, there should be $\tilde x_{i,j_tk_u}>0$, $s_j\times s_k$ positive entries in total. We make illustration for 2D case in \Cref{fig:gr 2d}. Note that the coordinates in transport are rearranged from the 2D coordinates in mesh. 
	
	\par Based upon the above arguments, we derive the GR subroutine for initialization; see \cref{alg:gr}. 
	
	\begin{algorithm}
		\caption{The GR initialization subroutine}
		\label{alg:gr}
		\begin{algorithmic}[1]
			\REQUIRE{Coarse mesh with $K$ elements $\{e_j\}$ and refined mesh with $\tilde K$ elements $\{\tilde e_k\}$; Solution of the previous step $(X_i)_{i=2}^N$; scaling factor $r>0$.}
			\FOR{$i=2,\ldots,N$}
			\FOR{$j=1,\ldots,K$}
			\FOR{$k=1,\ldots,K$}
			\IF{$x_{i,jk}>0$}
			\STATE{Find $\tilde e_{j_t},~t=1,\ldots,s_j$ such that $e_j=\bigcup_{t=1}^{s_j}\tilde e_{j_t}$.}
			\STATE{Find $\tilde e_{k_u},~u=1,\ldots,s_k$ such that $e_k=\bigcup_{u=1}^{s_k}\tilde e_{k_u}$.}
			\STATE{Set $\tilde x_{i,j_t,k_u}=r\cdot x_{i,jk}$ for $t\in\{1,\ldots,s_j\},~u\in\{1,\ldots,s_k\}$.}
			\ENDIF
			\ENDFOR
			\ENDFOR
			\ENDFOR
			\RETURN{$(\tilde X_i)_{i=2}^N\in\big(\R^{\tilde K\times\tilde K}\big)^{N-1}$, where for any $i$, $\tilde X_i=(\tilde x_{i,mn})_{mn}$.}
		\end{algorithmic}
	\end{algorithm}
	
	\subsection{Local Solver}\label{subsec:Algorithm}
	
	The global solver and the GR subroutine make brute-force globally solving large-scale \cref{eqn:penalty problem} unnecessary. Instead, we only need to provide a local solver (see \cref{local line} in \Cref{frame:ggr}). We assume the procedure is in the $(l+1)$-th iteration of \Cref{frame:ggr}. This is the same in the sequel whenever talking about local solver.
	
	\par Regarding algorithm design, the block structure of \cref{eqn:penalty problem} reminds us of using splitting type methods. One natural choice is an $(N-1)$-block cyclic \texttt{PBCD}; see \cref{frame:PBCD scheme}. In \texttt{PBCD}, the $i$-th block problem merely depends on the $i$-th block variable $X_i$, while keeping other block variables their latest values, where $f_{\beta^{(l+1)}}$ is defined in \cref{eqn:penalty problem}. Moreover, proximal term $\snorm{X_i-X_i^{(l,k)}}_\F^2$ is added to the objective function such that the block problem admits unique global solution. Here $\sigma>0$ is the proximal parameter.
	
	\begin{algorithm}
		\caption{\texttt{PBCD} for \cref{eqn:penalty problem}}
		\label{frame:PBCD scheme}
		\begin{algorithmic}[1]
			\REQUIRE{$R^{(l+1)},X_i^{(l,0)}\in\mathbb{R}^{K^{(l+1)}\times K^{(l+1)}},~i=2,\ldots,N$;\\\qquad\quad $e^{(l+1)},\varrho^{(l+1)}\in\mathbb{R}^{K^{(l+1)}};~\beta^{(l+1)},\sigma>0$.}
			\STATE{Set $k:=0$.}
			\WHILE{certain stopping criteria are not satisfied}
			\STATE{For $i=2,\ldots,N$, inexactly solve
				\begin{equation}\label{eqn:block problem} \min_{X_i\in\calS^{(l+1)}} f_{\beta^{(l+1)}}\big(X_{<i}^{(l,k+1)},X_i,X_{>i}^{(l,k)}\big)+\frac{\sigma}{2}\snorm{X_i-X_i^{(l,k)}}_{\F}^2\end{equation}
				to obtain $X_i^{(l,k+1)}$.}
			\STATE{$k:=k+1$.}
			\ENDWHILE
			\RETURN $(X_i^{(l+1)})_{i=2}^N:=(X_i^{(l,k)})_{i=2}^N\in\big(\R^{K^{(l+1)}\times K^{(l+1)}}\big)^{N-1}$.
		\end{algorithmic}
	\end{algorithm}
	
	\par Zooming in on \cref{eqn:block problem} in \cref{frame:PBCD scheme}, we find that the block problems are essentially strongly convex quadratic programmings, or more precisely, projecting a point onto $\calS^{(l+1)}$. Since the projection does not possess a closed-form expression, iterate infeasibility w.r.t. $\big(\calS^{(l+1)}\big)^{N-1}$ is inevitable in \cref{frame:PBCD scheme}. On the one hand, this brings difficulties in analyzing the convergence; see \cref{sec:Convergence Analysis}. On the other hand, there exist abundant algorithm resources for solving \eqref{eqn:block problem}. For instance, we can extend the semismooth Newton-CG method proposed in \cite{li2020efficient} to tackle \cref{eqn:block problem}; see more discussions in \cref{subsec:default solver}.
	
	\section{Convergence Analysis}\label{sec:Convergence Analysis}
	
	In this section, we show the convergence of the \texttt{PBCD} algorithm (\cref{frame:PBCD scheme}) to first-order stationary points (Karush-Kuhn-Tucker (KKT) points) or  global solutions of \cref{eqn:penalty problem} in different settings. The definition of KKT points for \cref{eqn:penalty problem} can be found in the supplementary material. 
	
	\par Since the convergence results are independent of the skeleton of the GGR approach (\Cref{frame:ggr}), we omit outer iteration index in the superscripts as well as the specification on the problem size; e.g., use $X_i^{(k)}$ instead of $X_i^{(l,k)}$. For the sake of brevity, we adopt the abbreviation $Z^{(k)}:=(X_i^{(k)})_{i=2}^N$ and $F(Z):=f_{\beta}(Z)+\sum_{i=2}^N\delta_{\calS}(X_i)$, where $f_{\beta}$ is defined in \cref{eqn:penalty problem}.
	
	\par When the block problems are exactly solved, we can directly follow the results in \cite{Xu2013A} and obtain the following theorem. 
	
	\begin{theorem}[Global Convergence of \cref{frame:PBCD scheme} -- Exact Version]\label{thm:global convergence exact}
		Let $\sigma>0$, and $\{Z^{(k)}\}$ be the sequence generated by \cref{frame:PBCD scheme} where block problems are exactly solved. Then $\{Z^{(k)}\}$ converges to a KKT point of \cref{eqn:penalty problem}. Moreover, $\{Z^{(k)}\}$ converges to a global minimizer of \cref{eqn:penalty problem} if the initial point $Z^{(0)}\in\calS^{N-1}$ is sufficiently close to some global minimizer. 
	\end{theorem}
	
	\par Since block exact solutions are not available in our case, we turn to study the global convergence property of \cref{frame:PBCD scheme} allowing block problems to be solved inexactly; in particular, the iterates are permitted to be infeasible w.r.t. $\calS^{N-1}$. 
	
	\par In the nonconvex context, existing convergence analyses of the \texttt{PBCD} algorithm restrict the iterates to be feasible,
	regardless of the complicate feasible set \cite{bolte2014proximal,gur2020convergent}. Limitations as they have, their analyses pave way for our study. Before presenting the results, we define the block optimal sequence $\{(\bar X_i^{(k)})_{i=2}^N\}$ as follows: for any $i\in\{2,\ldots,N\}$,
	\begin{equation}
		\bar X_i^{(k+1)}:=\argmin_{X_i\in\calS}f_{\beta}\big(X_{<i}^{(k+1)},X_i,X_{>i}^{(k)}\big)+\frac{\sigma}{2}\snorm{X_i-X_i^{(k)}}_\F^2.
		\label{eqn:optimal solution sequence}
	\end{equation}
	In other words, $\bar X_i^{(k+1)}$ is the unique global solution of the $i$-th block problems (\cref{eqn:block problem} in \cref{frame:PBCD scheme}). For any $k$, let $\bar Z^{(k)}:=(\bar X_i^{(k)})_{i=2}^N$. To facilitate analysis, we need assumptions on the local solver and energy value sequence.
	
	\begin{assumption}[Assumptions on the Local Solver and Energy Sequence]\label{assume:iterate}
		\mbox{}
		\begin{enumerate}[label=\textup{(}\arabic*\textup{)}]
			\item $\{F(\bar Z^{(k)})\}$ is non-increasing;\label{assume:decreasing energy}
			\item $\sum_{k=1}^\infty \Vert \bar Z^{(k)}-Z^{(k)}\Vert_\F<\infty$.\label{assume:local solver}
		\end{enumerate}
	\end{assumption}
	
	\par Since $F$ is continuous over the compact $\calS^{N-1}$, $F$ must attains its infimum in $\calS^{N-1}$. Hence \Cref{assume:iterate} (1) actually yields that the sequence $\{F(\bar Z^{(k)})\}$ converges to some $\bbar F\ge0$. \Cref{assume:iterate} (2) lays restrictions on the local solver, in that the block problems \cref{eqn:block problem} are solved more and more accurately.
	
	\par Since the analysis is rather involved, we give a rough statement of the convergence result for the \texttt{PBCD} algorithm below. The formal statement and proof are left to the supplementary material.
	\begin{theorem}[Convergence Property of \cref{frame:PBCD scheme} -- Inexact Version]\label{thm:rough convergence property}
		Suppose \Cref{assume:iterate} holds and $\sigma$ is sufficiently large. Let $\{Z^{(k)}\}$ be the sequence generated by \cref{frame:PBCD scheme}, $\{\bar Z^{(k)}\}$ be the  sequence defined in \cref{eqn:optimal solution sequence}. Assume also that $F(\bar Z^{(k)})>\bbar F,~\forall~ k\ge1$. Then 
		\begin{enumerate}[label=\textup{(}\arabic*\textup{)}]
			\item the sequence $\{Z^{(k)}\}$ converges to a KKT point of \cref{eqn:penalty problem};
			\item if further $\sum_{k=1}^\infty\Vert \bar Z^{(k)}-Z^{(k)}\Vert_\F$ is small enough, $Z^{(0)}$ is feasible and sufficiently close to some global minimizer, then $\{Z^{(k)}\}$ converges to a global minimizer of \cref{eqn:penalty problem}. 
		\end{enumerate}
	\end{theorem}
	
	\section{Numerical Experiments}\label{sec:Numerical Experiments}
	
	In this section, we validate the proposed GGR approach via numerical simulations on several typical systems, including both 1D and 2D systems. During the experiments, we mainly monitor the repulsive energy $f$ in \cref{eqn:complicate MPCC} (denoted by E) . We also calculate the approximated transport maps $\{T_i^K\}_{i=2}^N$ as in \cref{eqn:approx transport maps}, and in turn evaluate the quality of solution through the average error (denoted by err)
	$$\mathrm{err}(K,\Omega):=\frac{1}{K\abs{\Omega}}\sum_{i=1}^K\sum_{j=2}^N\abs{T_j(\aaa_i)-T_j^K(\aaa_i)},$$
	if the optimal transport maps $\{T_i\}_{i=2}^N$ \cref{eqn:Monge formulation} are already available. We refer interested readers to supplementary material for the numerical comparison among local solvers proposed in \cite{byrd2006knitro,fletcher2004solving,waltz2006interior} and \texttt{PBCD}.
	
	\par All the numerical experiments presented here are run in a platform with Intel(R) Xeon(R) Gold 6242R CPU @ 3.10GHz and 510GB RAM running MATLAB R2018b under Ubuntu 20.04.
	
	\subsection{Default Settings}\label{subsec:default solver}
	\mbox{}
	
	\par\textbf{Global solver}. Considering the applicability and efficiency, we take the stochastic method, \textit{random multi-start}, as global solver for 1D systems, and employ software \texttt{BARON} for 2D systems. The implementation of random multi-start follows \cite{hickernell1997simple}. The software \texttt{BARON} invokes efficient random multi-start procedures initially, and then carries out \textit{branch-and-bound} and \textit{cutting plane} algorithm for global optimization. Version 21.1.13 of \texttt{BARON} is available in the downloadable AMPL system \cite{fourer2003ampl}.
	
	\vskip 0.2cm
	
	\par\textbf{Details in \texttt{PBCD}}. We adapt the semi-smooth Newton-CG (\texttt{SSNCG}) method in \cite{li2020efficient} to solve the dual block problems. A general iteration in \texttt{SSNCG} consists of approximately solving a sparse symmetric positive definite linear system of the form 
	$$\big(\mathcal{V}^{(l+1)}+\eps I\big)d+r^{(l+1)}=0,\quad d\in\mathrm{Im}\big(\mathcal{B}^{(l+1)}\big),$$
	and then performing line searches along direction $d$ for a sufficient reduction on dual objective. Here $\mathcal{B}^{(l+1)}$ is a linear operator defined as
	\begin{equation}\label{eqn:linear operator B}\mathcal{B}^{(l+1)}(W)=[(e^{(l+1)})^\T W^\T\quad(\varrho^{(l+1)})^\T\Xi^{(l+1)} W\quad\trace(W)],\end{equation}
	$\mathcal{V}^{(l+1)}$ is a positive semidefinite operator associated with $\mathcal{B}^{(l+1)}$, $\eps>0$, $I$ stands for an identity matrix in proper dimension for convenience, and $r^{(l+1)}$ is the residual vector. In our context, the linear system can be solved quickly to desirability by the preconditioned conjugate gradient method equipped with block Jacobi preconditioner. 
	
	\vskip 0.2cm
	
	\par\textbf{Parameter setting}. In the GR subroutine, we set scaling value $r=1$.
	In all experiments, we fix $\sigma=10^{-3},~\mathrm{maxit}=10^6$ in \texttt{PBCD}. For different $K$, we choose $\beta$ according to \Cref{tab:beta}. 
	\begin{table}[htbp]
		\small
		\centering
		\caption{The value of $\beta$ for different $K$}
		\label{tab:beta}
		\begin{tabular}{cccccc}
			\toprule
			$K$ & $(0,10)$ & $[10,36)$ & $[36,80)$ & $[80,160)$ & $[160,320)$ \\\midrule
			$\beta$ & $2^2$ & $2^{1}$ & $2^{0}$ & $2^{-2}$ & $2^{-3}$ \\\midrule
			$K$ & $[320,640)$ & $[640,1280)$ & $[1280,2560)$ & $[2560,5120)$ & $[5120,\infty)$ \\\midrule
			$\beta$ & $2^{-4}$ & $2^{-5}$ & $2^{-6}$ & $2^{-7}$ & $2^{-8}$ \\
			\bottomrule
		\end{tabular}
	\end{table}
	We invoke \texttt{SSNCG} for block problems in \texttt{PBCD}. The maximum \texttt{SSNCG} iteration number is set to $\mathrm{maxitSSN}=10^5$. We start \texttt{SSNCG} from zero point in the first call; after that, we perform \textit{wart start}. 
	
	\vskip 0.2cm
	
	\par\textbf{Stopping criteria}. We terminate \texttt{SSNCG} if feasibility violation
	\begin{equation}
		\mathrm{Feas}(Z)=\sum_{i=2}^N\big\Vert\mathcal{B}^{(l+1)}(X_i)-b^{(l+1)}\big\Vert_2
		\label{eqn:primal infeas}
	\end{equation}
	is smaller than $\eps_{\mathrm{inner}}=10^{-9}$ in all cases, where $\mathcal{B}^{(l+1)}$ is the linear operator defined in \cref{eqn:linear operator B} and $b^{(l+1)}=[\one^\T~(\varrho^{(l+1)})^\T~0]^\T$. We stop \texttt{PBCD} when the scaled difference of two consecutive iterate $\sqrt{\sigma}\snorm{Z^{(l,k+1)}-Z^{(l,k)}}_\F$ is less than a prescribed value $\eps_{\mathrm{outer}}$, which is chosen as
	\begin{equation}
		\label{eqn:outer kkt threshold}
		\eps_{\mathrm{outer}}=\left\{\begin{array}{ll}
			10^{-8}, & K^{(l+1)}\in(0,200],\\
			10^{-6}, & K^{(l+1)}\in(200,2000],\\
			10^{-5}, & K^{(l+1)}\in(2000,10000],\\
			10^{-4}, & K^{(l+1)}\in(10000,\infty),
		\end{array}\right.
	\end{equation}
	or, when the absolute value of the difference between two consecutive energies is less than $10^{-8}$.

	\subsection{Numerical Results on 1D Systems}
	\label{subsec:1D results}
	
	We first consider some typical 1D systems with our GGR approach. 
	In the simulations, we use ``equal-mass" discretization of the marginal for the initial mesh, in that each element in mesh carries the same marginal mass. This can be achieved cheaply and exactly for 1D systems. The meshes are refined uniformly afterwards. 
	
\begin{figure}[htbp]
		\centering
		\includegraphics[width=.8\linewidth]{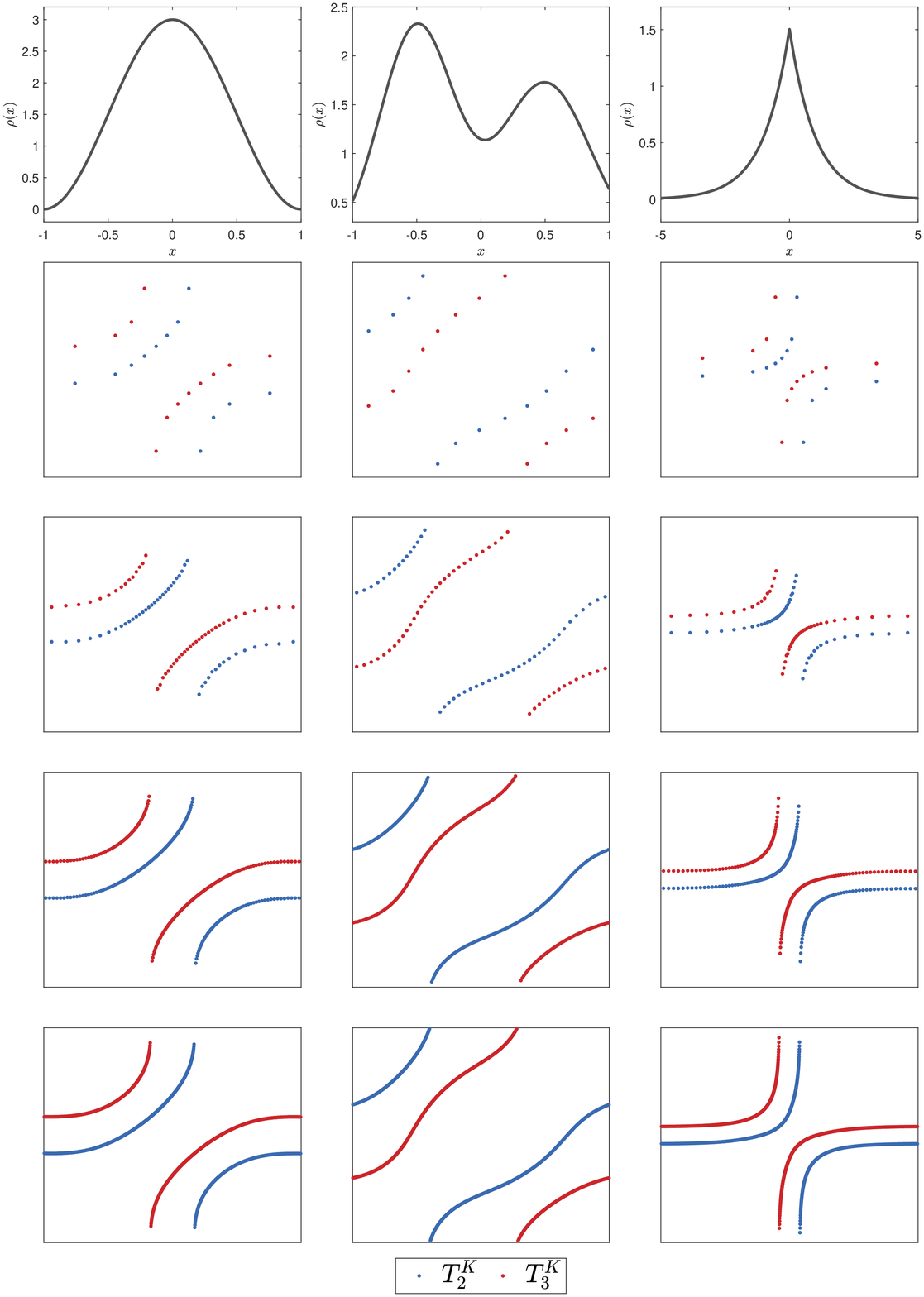}
		\caption{Marginals (the first row) and approximated transport maps (the remaining four rows) in 1D $N=3$ systems; left to right: systems with $\rho_1,\rho_2,\rho_3$. The maps in the last four rows correspond to rows in \Cref{tab:1D} \textup{(a)} where $K=12,48,192,768$. The blue and red dots are images of $T_2^K$ and $T_3^K$ over grid barycenters, respectively.}
		\label{fig:1D N3}
	\end{figure}

	The first three systems under consideration all consist of 3 particles ($N=3$), whose single-electron densities (marginals) are given by
	\begin{align*}
		& \rho_1(x) = c_1 \big(\cos(\pi x)+1\big), \quad & \Omega = [-1,1],    
		\\
		& \rho_2(x) = c_2 \big(2e^{-6(x+0.5)^2}+1.5e^{-4(x-0.5)^2}\big), \quad & \Omega = [-1,1],   
		\\
		& \rho_3(x) = c_3 e^{-\abs{x}}, \quad & \Omega = [-5,5],
	\end{align*}
	respectively, with $c_i,~i=1,2,3,$ the normalization constants such that $\int_{\Omega}\rho_i(x)\dd x = 3$.
	The number of grid points used for the initial meshes is $K^{(0)}=12$ for all three systems. The single-electron densities (marginals) and the approximated transport maps $\{T_i^K\}$ \cref{eqn:approx transport maps} are shown in \Cref{fig:1D N3}.
	The convergence of the GGR approach can be observed as the meshes being refined. Note that explicit solutions of the original MMOT problems are known for 1D systems \cite{seidl07strictly}, our results can match the theory perfectly. We also list the output energies and the calculated average errors (the ``err\_e'' column) at each step in \Cref{tab:1D} (a), supporting the efficiency of our approach.

	The second set includes three systems, each of which contains 7 particles ($N=7$). 
	Note that this particle number is already intractable if one tries to solve the original MMOT problem \cref{MMOT:Coulomb} directly.
	The single-electron densities (marginals) are given by
	\begin{align*}
		& \rho_4(x) = c_4  e^{-x^2/\sqrt{\pi}}, \quad & \Omega = [-2,2],    
		\\
		& \rho_5(x) = c_5 \big(e^{-3(x+3)^2}+e^{-3(x+2)^2}+e^{-2(x+1)^2}+e^{-x^2}+e^{-2(x-1)^2} &
		\\ & \qquad\qquad +e^{-3(x-2)^2}+e^{-3(x-3)^2}\big),  & \Omega = [-4,4], 
		\\
		& \rho_6(x) = c_6 \big(e^{-8(x+2.7)^2}+e^{-8(x+2.025)^2}+e^{-8(x+1.35)^2}+e^{-8(x+0.675)^2} & 
		\\
		& \qquad\qquad +e^{-5(x-0.5)^2}+e^{-5(x-1.5)^2}+e^{-5(x-2.5)^2}\big), \quad & \Omega = [-3,3],
	\end{align*}
	respectively, with $c_i,~i=4,5,6,$ the normalization constants such that $\int_{\Omega}\rho_i(x)\dd x = 7$.
	The last two examples can be viewed as systems with localized electrons, where each Gaussian represents the distribution of an electron.
	The number of grid points used for the initial meshes is $K^{(0)}=14$ for these three systems.
	
\begin{figure}[htbp]
		\centering
		\includegraphics[width=.8\linewidth]{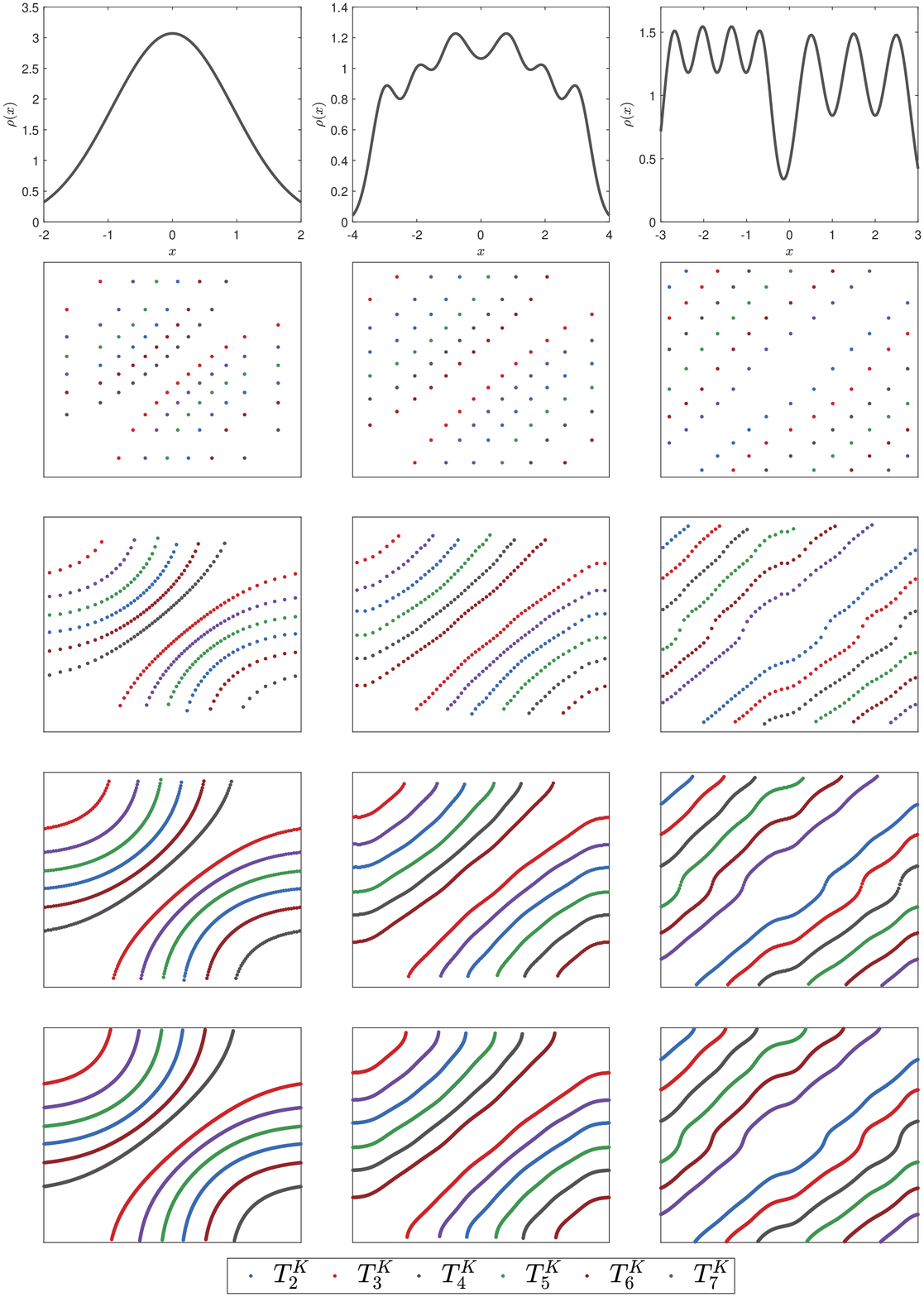}
		\caption{Marginals (the first row) and approximated transport maps (the remaining four rows) in 1D $N=7$ systems; left to right: systems with $\rho_4,\rho_5,\rho_6$. The maps in the last four rows correspond to rows in \Cref{tab:1D} \textup{(b)} where $K=14,56,224,896$. The blue, red, black, green, brown and purple dots are images of $T_i^K,~i=2,\ldots,7$ over grid barycenters, respectively.}
		\label{fig:1D N7}
	\end{figure}

	We show the single-electron densities (marginals) and the convergence of the GGR approach during the mesh refinements in \Cref{fig:1D N7}.
	We also list the output energies and the calculated average errors at each step in \Cref{tab:1D} (b).
	We observe from the numerical results that the iterates given by our GGR approach converge well to the correct solutions.
	\begin{table}[htbp]
		\centering
		\caption{Output energies and calculated average errors of the GGR approach on 1D systems}
		\label{tab:1D}
		\resizebox{\textwidth}{16mm}{
			\subfloat[$N=3$]{
				\begin{tabular}{l|rrrr||rrrr||rrrr}
					\toprule
					\multirow{2}{*}{Step} & \multicolumn{4}{c||}{System 1} & \multicolumn{4}{c||}{System 2} & \multicolumn{4}{c}{System 3} \\\cmidrule{2-13}
					& $K$ & E & err\_s & err\_e & $K$ & E & err\_s & err\_e & $K$ & E & err\_s & err\_e \\\midrule
					GGR\_Init  & 12 & 18.114 & - & 0.031 & 12 &  12.211 & - & 0.012 & 12 &  6.024 & - & 0.040 \\
					GGR\_LS(1) & 24 & 18.911 & 0.049 & 0.013 & 24 & 12.373 & 0.041 & 0.011 & 24 & 6.318 & 0.053 & 0.018 \\
					GGR\_LS(2) & 48 & 19.004 & 0.022 & 0.009 & 48 & 12.367 & 0.026 & 0.012 & 48 & 6.389 & 0.027 & 0.013 \\
					GGR\_LS(3) & 96 & 19.019 & 0.014 & 0.004 & 96 & 12.360 & 0.017 & 0.009 & 96 & 6.400 & 0.026 & 0.011 \\
					GGR\_LS(4) & 192 & 19.021 & 0.007 & 0.003 & 192 & 12.358 & 0.010 & 0.003 & 192 & 6.403 & 0.013 & 0.001 \\
					GGR\_LS(5) & 384 & 19.022 & 0.007 & 0.002 & 384 & 12.358 & 0.005 & 0.001 & 384 & 6.404 & 0.003 & 0.000 \\
					GGR\_LS(6) & 768 & 19.022 & 0.004 & 0.001 & 768 & 12.357 & 0.003 & 0.001 & 768 & 6.404 & 0.001 & 0.000 \\
					\bottomrule
		\end{tabular}}}
		
		\resizebox{\textwidth}{16mm}{
			\subfloat[$N=7$]{
				\begin{tabular}{l|rrrr||rrrr||rrrr}
					\toprule
					\multirow{2}{*}{Step} & \multicolumn{4}{c||}{System 4} & \multicolumn{4}{c||}{System 5} & \multicolumn{4}{c}{System 6} \\\cmidrule{2-13}
					& $K$ & E & err\_s & err\_e & $K$ & E & err\_s & err\_e & $K$ & E & err\_s & err\_e \\\midrule
					GGR\_Init  & 14 & 189.626 & - & 0.018 & 14 & 80.266 & - & 0.021 & 14 & 91.536 & - & 0.016 \\
					GGR\_LS(1) & 28 & 193.703 & 0.027 & 0.023 & 28 & 82.199 & 0.024 & 0.012 & 28 & 93.056 & 0.022 & 0.025 \\
					GGR\_LS(2) & 56 & 193.312 & 0.019 & 0.026 & 56 & 81.937 & 0.012 & 0.012 & 56 & 92.458 & 0.025 & 0.019 \\
					GGR\_LS(3) & 112 & 193.128 & 0.022 & 0.016 & 112 & 81.854 & 0.010 & 0.014 & 112 & 92.245 & 0.019 & 0.013 \\
					GGR\_LS(4) & 224 & 193.066 & 0.015 & 0.010 & 224 & 81.817 & 0.014 & 0.011 & 224 & 92.185 & 0.020 & 0.007 \\
					GGR\_LS(5) & 448 & 193.044 & 0.007 & 0.004 & 448 & 81.808 & 0.009 & 0.002 & 448 & 92.171 & 0.007 & 0.003 \\
					GGR\_LS(6) & 896 & 193.039 & 0.002 & 0.002 & 896 & 81.806 & 0.001 & 0.002 & 896 & 92.167 & 0.003 & 0.001 \\
					\bottomrule
		\end{tabular}}}
	\end{table}
	
	\par To show that the GR subroutine yields high-quality initialization, we compute the average errors of the initial points (the ``err\_s'' column) as well; the notation ``-'' in the GGR\_Init step indicates no initial point is fed to global solver. The decreasing err\_s's underline the efficacy of the GR subroutine, which boosts the GGR approach and helps us find global solutions. Incidentally, the comparison between err\_s and err\_e in the same row highlights the improvements due to the local solver \texttt{PBCD}. Meanwhile, one can see that err\_e is sometimes slightly larger than err\_s. In these cases, \texttt{PBCD} eliminates infeasibility while inheriting the high quality of initial points.

	\subsection{Numerical results on 2D systems}\label{subsec:2D results}
	
	We then consider some 2D systems with the GGR approach. 
	We use the finite elements package \texttt{FreeFEM} \cite{freefem} to generate the initial meshes for the marginal discretization.
	The meshes are non-uniform such that every element carries almost the same mass.
	In the later loops of the GGR approach, each element is refined in the same manner.
	
	\par
\begin{figure}[htbp]
		\centering
		\includegraphics[width=.9\linewidth]{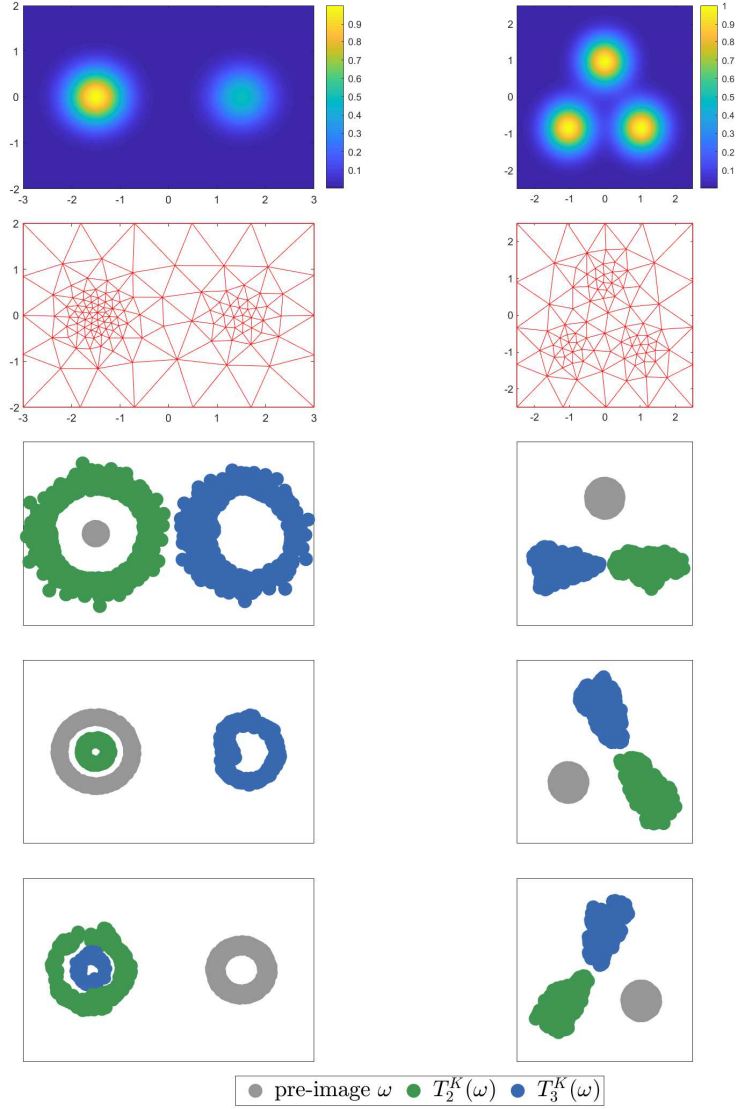}
		\caption{Contours of marginals (the first row), initial meshes (the second row) and slices of approximated transport maps (the third-fifth row) in 2D systems; left to right: systems with $\rho_7,\rho_8$. In system 7-8, we calculated $K$ to 15360 and 10880, respectively. The gray, blue and green circles are pre-images $\omega\subseteq\Omega$, $T_2^K(\omega)$ and $T_3^K(\omega)$, respectively.}
		\label{fig:2d N3}
	\end{figure}

 The two systems under consideration both consist of 3 particles ($N=3$), whose single-electron densities (marginals) are given by
	\begin{align*}
		& \rho_7(x,y) = c_7\big(e^{-2.5\abs{(x,y)-(-1.5,0)}^2}+0.5e^{-2.5\abs{(x,y)-(1.5,0)}^2}\big), \quad & \Omega = [-3,3]\times[-2,2],\\
		& \rho_8(x,y) = c_8\big(e^{-2.5\abs{(x,y)-(-1.032,-0.84)}^2)}+e^{-2.5\abs{(x,y)-(0,0.96)}^2} &\\
		& \qquad\qquad +e^{-2.5\abs{(x,y)-(1.032,-0.84)}^2}\big), \quad & \Omega = [-2.5,2.5]^2,
	\end{align*}
	respectively, with $c_i,~i=7,8,$ the normalizing factors such that $\int_\Omega\rho_i(x,y)\dd x\dd y=7$. For the first 2D system, $\rho_7$ corresponds to a system that has two electrons located on the left part of $\Omega$ (represented by the first Gaussian centered at $(-1.5,0)$), and the third electron located on the right part (represented by the second Gaussian centered at $(1.5,0)$).
	For the second 2D system, $\rho_8$ corresponds to a system that has three electrons concentrated on three different sites $(-1.032,-0.84)$, $(0,0.96)$ and $(1.032,-0.84)$ (represented by three Gaussians), respectively. The electron densities (marginals) and the corresponding initial meshes (obtained by \texttt{FreeFEM}) are shown in the first two rows of \Cref{fig:2d N3}.
	The numbers of grid points used for the initial meshes are $K^{(0)}=240$ for $\rho_7$ and $K^{(0)}=170$ for $\rho_8$, respectively. After three steps in \Cref{frame:ggr}, we reach $K^{(3)}=15360$ for $\rho_7$ and $K^{(3)}=10880$ for $\rho_8$, respectively.
	
	\begin{table}[htbp]
		\footnotesize
		\centering
		\caption{Output energies of the GGR approach on 2D systems}
		\label{tab:2D}
		\begin{tabular}{l|rr||rr}
			\toprule
			\multirow{2}{*}{Step} & \multicolumn{2}{c||}{System 7} & \multicolumn{2}{c}{System 8}\\\cmidrule{2-5}
			& $K$ & E & $K$ & E \\\midrule
			GGR\_Init & 240 & 9.503 & 170 & 9.491 \\
			GGR\_LS(1) & 960 & 9.577 & 680 & 9.533 \\
			GGR\_LS(2) & 3840 & 9.598 & 2720 & 9.543 \\
			GGR\_LS(3) & 15360 & 9.604 & 10880 & 9.546 \\
			\bottomrule
		\end{tabular}
	\end{table}
	
	\par 
	To show our results on the 2D systems, in the remaining three rows of \Cref{fig:2d N3}, we plot the images of the \textit{barycenters} of triangular elements within some given regions $\omega\subseteq\Omega$, where $T^K_2$ and $T^K_3$ are the approximated transport maps \cref{eqn:approx transport maps} given by the GGR approach. For the two-Gaussian system $\rho_7$, the pictures show that: if the first electron is around the left Gaussian center, then the third electron will go to the region near the right Gaussian center, and the second electron will lie in the left part (to satisfy the marginal constraints) but stay away from the first one ($\omega$ and $T^K_2(\omega)$ lie in two different regions around the left Gaussian center); if one electron is located around the right Gaussian center, then the other two electrons will be around the left Gaussian center while keeping distance away from each other.
	For the three-Gaussian system $\rho_8$, we can see that if one electron is located around one of the Gaussian centers, then the other two electrons go to the other two Gaussian centers, respectively. We also list the output energies at each step in \Cref{tab:2D}. Though there are no theoretical results for comparison, our simulations match physical intuitions quite well and can support the reliability of our approach. 
	
	\section{Conclusions}\label{sec:conclusion}
	In the present work, we consider the MMOT problem with Coulomb cost arising in quantum physics. The Monge-like ansatz tides us over curse of dimensionality, in that the number of unknowns scales linearly w.r.t. the number of electrons, however resulting in MPGCC. In quest for global solutions, we propose a global optimization approach GGR for dealing with the derived MPGCC. The GGR approach solves the problem step by step along with the process of mesh refinement, and is equipped with an initialization subroutine such that global solutions are amenable to the proposed local solver \texttt{PBCD}. The convergence property of \texttt{PBCD} is established in the presence of iterate infeasibility. We corroborate the merits of the GGR approach with numerical simulations on several typical 1D and 2D physical systems. Notably, we obtain solutions with high resolution in the 1D cases, and visualize the optimal transport maps in the 2D context.
	
	\appendix
	
	\section{Discretization of (\ref{eqn:optimization formulation})}\label{appsec:discretization and error analysis}
	
	For the repulsive energy in \cref{eqn:optimization formulation}, we have for any $i\in\{2,\ldots,N\}$,
	$$\int_{\Omega}\int_\Omega\frac{\rho(\rr)\gamma_i(\rr,\rr')}{\abs{\rr-\rr'}}\dd\rr\dd\rr'=\sum_{j,k}\int_{e_k}\int_{e_j}\frac{\rho(\rr)\gamma_i(\rr,\rr')}{\abs{\rr-\rr'}}\dd\rr\dd\rr'.$$
	Note when $k=j$, the integral explodes and hence we impose $x_{kk}^i=0,~\forall~ k$ as extra constraints to avoid numerical instability. In the sequent derivation, we take $\gamma_i(\rr,\rr')=0$ whenever $\rr$ and $\rr'$ belong to the same element:
	\begin{align}
		\int_{\Omega}\int_\Omega\frac{\rho(\rr)\gamma_i(\rr,\rr')}{\abs{\rr-\rr'}}\dd\rr\dd\rr'=&\sum_{j\ne k}\int_{e_k}\int_{e_j}\frac{\rho(\rr)\gamma_i(\rr,\rr')}{\abs{\rr-\rr'}}\dd\rr\dd\rr'\nonumber\\
		=&\sum_{j\ne k}\varrho_jx_{i,jk}\int_{e_k}\int_{e_j}\frac{1}{\abs{\rr-\rr'}}\dd\rr\dd\rr'+O(h)\nonumber\\
		=&\sum_{j\ne k}\varrho_jx_{i,jk}r_{jk}\abs{e_j}\abs{e_k}+O(h)=\inner{X_i,\Lambda\Xi C\Xi}+O(h)\label{eqn:appendix linear term},
	\end{align}
	where $h:=\snorm{e}_\infty$ represents the size of the largest element. 
	By similar reasoning, we can write for any $i,j\in\{2,\ldots,N\}:i\ne j$,
	\begin{align}
		&\int_{\Omega}\int_\Omega\int_\Omega\frac{\rho(\rr)\gamma_i(\rr,\rr')\gamma_j(\rr,\rr'')}{\abs{\rr'-\rr''}}\dd\rr\dd\rr'\dd\rr''\nonumber\\
		=&\sum_{m,n,t:n\ne t}\varrho_mx_{i,mn}x_{j,mt}r_{nt}\abs{e_m}\abs{e_n}\abs{e_t}+O(h)\label{eqn:appendix quadratic term}\\
		=&\inner{X_i,\Xi\Lambda X_j\Xi C\Xi}+O(h).\nonumber
	\end{align}
	Note that we have excluded $n=t$ cases and impose $\inner{X_i,X_j}=0$ as extra complementarity constraints. By \cref{eqn:appendix linear term} and \cref{eqn:appendix quadratic term}, the repulsive energy in \cref{eqn:optimization formulation} can be approximated by
	$$\sum_{2\le i\le N}\inner{X_i,\Lambda\Xi C\Xi}+\sum_{i<j}\inner{X_i,\Xi\Lambda X_j\Xi C\Xi},$$
	with error depending on the size of the largest element.
	\par Regarding normalizing and marginal constraints in \cref{conditions:gammas}, we can see from similar derivation that, for any $i\in\{2,\ldots,N\}$,
	$$\begin{aligned}
		1=\frac{1}{\abs{e_j}}\int_{e_j}1\dd\rr&=\frac{1}{\abs{e_j}}\int_{e_j}\int_{\Omega}\gamma_i(\rr,\rr^\prime)\dd\rr'\dd\rr=\sum_{k=1}^Kx_{i,jk}\abs{e_k},~\forall j,\\
		\varrho_k=\frac{1}{\abs{e_k}}\int_{e_k}\rho(\rr')\dd\rr'&=\frac{1}{\abs{e_k}}\int_{e_k}\int_{\Omega}\rho(\rr)\gamma_i(\rr,\rr')\dd\rr\dd\rr'=\sum_{j=1}^K\varrho_jx_{i,jk}\abs{e_j}+O(h),~\forall k.
	\end{aligned}$$
	Consequently, the constraints in \cref{conditions:gammas} can be approximated using
	$$X_ie=\one,~X_i^\T\Xi\rho=\rho,~\forall~ i\in\{2,\ldots,N\}.$$
	
	\normalem
	\bibliographystyle{siamplain}
	\bibliography{note}
	
\end{document}

\maketitle
	
	\par We provide technical results on \textup{\texttt{PBCD}} in \cref{sec:preliminaries,appsec:lemmas,appsec:detail proof}. \Cref{sec:preliminaries} is dedicated to some preliminaries, containing optimality conditions of the problem to be solved, some necessary definitions and key tools for convergence analysis. The proof ingredients and main results are detailed in \cref{appsec:lemmas,appsec:detail proof}, respectively.
	
	\par \Cref{subsec:comparison NLP} is devoted to the numerical comparison among local solvers, \texttt{Knitro}, \texttt{filtermpec} and our proposed \textup{\texttt{PBCD}}. The former two are nonlinear programming based solvers whose performance comparison on a benchmark is available online\footnote{See \url{http://plato.asu.edu/ftp/mpec.html}.}.
	
	\par The analysis and description hereinafter are irrelevant to the skeleton of the GGR approach. Therefore we make abbreviations on the superscripts and omit specification on size $K$. For example, we replace $Z^{(l,k)}$ with $Z^{(k)}$.
	
	\par For better readability, we restate some frequently referenced things. The problem to be solved by local solver \textup{\texttt{PBCD}} is 
	\begin{equation}
		\label{eqn:penalty problem restate}
		\min_{X_2,\ldots,X_N}~f_\beta(X_2,\ldots,X_N)\quad\st~\mathcal{B}(X_i)=b,~X_i\ge0,~i=2,\ldots,N,
	\end{equation}
	where the repulsive energy
	$$f_\beta(X_2,\ldots,X_N)=\sum_{2\le i\le N}\inner{X_i,\Lambda\Xi C\Xi}+\sum_{i<j}\inner{X_i,\Xi\Lambda X_j\Xi C\Xi}+\beta\sum_{i<j}\inner{X_i,X_j},$$
	the linear operator $\mathcal{B}:\R^{K\times K}\to\R^{2K+1}$ is defined as
	$$\mathcal{B}(W)=\big[e^\T W^\T\quad\varrho^\T\Xi W\quad\trace(W)\big]^\T,\quad\forall~ W\in\R^{K\times K},$$
	and the vector $b=[\one^\T\quad\varrho^\T\quad 0]^\T\in\R^{2K+1}$. 
	
	\par The \textup{\texttt{PBCD}} algorithm solves \cref{eqn:penalty problem restate} iteratively in a Gauss-Seidel style, i.e. in each cycle inexactly solves the block problem
	\begin{equation}\label{eqn:PBCD scheme}\min_{X_i\in\calS}f_\beta(X_{<i}^{(k+1)},X_i,X_{>i}^{(k)})+\frac{\sigma}{2}\snorm{X_i-X_i^{(k)}}_\F^2\end{equation}
	to obtain $X_i^{(k+1)}$ for $i=2,\ldots,N$, where $\calS:=\{W\in\R^{K\times K}:\mathcal{B}(W)=b,~W\ge0\}$. The block optimal sequence $\{\bar Z^{(k)}:=(\bar X_i^{(k)})_{i=2}^N\}$ is cast as
	\begin{equation}
		\label{eqn:optimal solution sequence restate}
		\bar X_i^{(k+1)}:=\argmin_{X_i\in\calS}f_\beta(X_{<i}^{(k+1)},X_i,X_{>i}^{(k)})+\frac{\sigma}{2}\snorm{X_i-X_i^{(k)}}_\F^2,~\forall~ i=2,\ldots,N.
	\end{equation}
	
	\par To facilitate the convergence analysis on \textup{\texttt{PBCD}}, we make assumptions on energy value and local solver below.
	\begin{assumption}[Assumptions on the Local Solver and Energy Sequence]\label{assume:iterate restate}
		\begin{enumerate}[label=\textup{(}\arabic*\textup{)}]
			\item $\{F(\bar Z^{(k)})\}$ is non-increasing and hence converges to some $\bbar F>0$\label{assume:decreasing func 2};
			\item $\sum_{k=1}^\infty \Vert \bar Z^{(k)}-Z^{(k)}\Vert_\F<\infty$.\label{assume:local solver 2}
		\end{enumerate}
	\end{assumption}
	
	\section{Preliminaries}\label{sec:preliminaries}
	
	In what follows, we briefly go through the first-order stationary conditions of \cref{eqn:penalty problem restate} (\Cref{lem:KKT}) together with some notions and tools concerning the convergence analysis (e.g. \Cref{lem:KL}).
	
	\par Since in \cref{eqn:penalty problem restate}, the feasible set $\calS^{N-1}$ is a polyhedron, we naturally have the necessity of Karush-Kuhn-Tucker (KKT) conditions.
	
	\begin{lemma}[KKT Conditions of \cref{eqn:penalty problem restate}]\label{lem:KKT}
		If $(X_i)_{i=2}^N$ is a local minimizer of \cref{eqn:penalty problem restate}, there exist multipliers $\lambda_i\in\R^{2K+1},\Phi_i\in\R^{K\times K},~i=2,\ldots,N$ such that for any $i\in\{2,\ldots,N\}$,
		\begin{align*}
			&\Lambda\Xi C\Xi+\Xi\Lambda\lrbracket{\sum_{j\ne i}X_j}\Xi C\Xi+\beta\lrbracket{\sum_{j\ne i}X_j}-\mathcal{B}^*(\lambda_i)-\Phi_i=0,\\
			&\mathcal{B}(X_i)=b,~X_i,\Phi_i\ge0,~\Phi_i\circ X_i=0,
		\end{align*}
		where ``$\circ$'' denotes the Kronecker product between two matrices. The tuple $(X_i)_{i=2}^N$ satisfying the above relations is called a KKT point of \cref{eqn:penalty problem restate}.
	\end{lemma}
	
	The convergence analysis involves some definitions about subdifferentials and the well-known Kurdyka-\L ojasiewicz property. 
	
	\begin{definition}[Subdifferentials]\label{def:subdifferential}
		Let $G:V\to(-\infty,\infty]$ be a proper and lower semi-continuous function, where $V$ is an Euclidean space with $\inner{\cdot,\cdot}_V$ as its inner product and $\snorm{\cdot}_V=\sqrt{\langle\cdot,\cdot\rangle_V}$ as its norm.
		\begin{enumerate}[label=\textup{(}\arabic*\textup{)}]
			\item For a given $x\in\dom(G):=\{y:G(y)<\infty\}$, the Fr\'echet subdifferential of $G$ at $x$, denoted by $\widehat\partial G(x)$, is the set consisting of all vectors $u\in V$ satisfying
			$$\liminf_{y\ne x,y\to x}\frac{G(y)-G(x)-\inner{u,y-x}_V}{\snorm{y-x}_V}\ge0.$$
			When $x\notin\dom(G)$, we simply set $\widehat\partial G(x)=\emptyset$.
			\item The limiting-sudifferential of $G$ at $x\in V$, denoted by $\partial G(x)$, is defined as
			$$\partial G(x):=\lrbrace{u\in V:\exists~x^{(k)}\to x,~G(x^{(k)})\to G(x),u^{(k)}\in\widehat\partial G(x^{(k)})\to u}.$$
		\end{enumerate}
	\end{definition}
	
	\begin{remark}\label{rem:subdiff}
		\begin{enumerate}[label=\textup{(}\arabic*\textup{)}]
			\item Let $\{(x^{(k)},u^{(k)})\}\subseteq\mathrm{graph}(\partial G)$ be a sequence converging to $(x,u)$. Here $\mathrm{graph}(\partial G)$ is the graph of set-valued map $\partial G$ defined as $\mathrm{graph}(\partial G):=\{(x,u):u\in\partial G(x)\}$. By the definition of $\partial G(x)$, if $G(x^{(k)})$ converges to $G(x)$, then $(x,u)\in\mathrm{graph}(\partial G)$\label{rem:closedness of graph}.
			\item If $x\in V$ is a local minimizer of $G$, then $0\in\partial G(x)$; the point $x$ satisfying $0\in\partial G(x)$ is called a critical point of $G$.
			\item For the extended-value function
			$$F(X_2,\ldots,X_N):=f_\beta(X_2,\ldots,X_N)+\sum_{i=2}^N\delta_{\calS}(X_i),$$
			a tuple $(X_i)_{i=2}^N$ is a critical point of $F$ if and only if it is a KKT point of \cref{eqn:penalty problem restate}; see \cite[Proposition 6.6]{penot2012calculus2}.\label{rem:critical and KKT} 
		\end{enumerate}
	\end{remark}
	
	\begin{definition}[Function class $\Phi_\eta$]\label{def:function class}
		Let $\eta\in(0,\infty]$. $\Phi_\eta$ is the class of all concave and continuous functions $\varphi:[0,\eta)\to\mathbb{R}_+$ satisfying \textup{(i)} $\varphi(0)=0$; \textup{(ii)} $\varphi\in C^1((0,\eta))$ and is continuous at $0$; \textup{(iii)} $\varphi^\prime(s)>0,~\forall~ s\in(0,\eta)$.
	\end{definition}
	
	\par With the definition of limiting-subdifferential (\Cref{def:subdifferential}) and function class $\Phi_\eta$ (\Cref{def:function class}), we can state the Kurdyka-\L ojasiewicz (K\L) property as follows.
	
	\begin{definition}[K\L~property]\label{def:KL property}
		Let $G:V\to(-\infty,\infty]$ be proper and lower semi-continuous, where $V$ is an Euclidean space with $\inner{\cdot,\cdot}_V$ as its inner product and $\snorm{\cdot}_V=\sqrt{\langle\cdot,\cdot\rangle_V}$ as its norm.
		\begin{enumerate}[label=\textup{(}\arabic*\textup{)}]
			\item $G$ is said to have the K\L~property at $\bar x\in\dom(\partial G):=\lrbrace{x\in V:\partial G(x)\ne\emptyset}$ if there exist a scalar $\eta\in(0,\infty]$, a neighborhood $\mathcal{U}$ of $\bar x$ and $\varphi\in\Phi_\eta$, such that for all $x\in \mathcal{U}\bigcap\{y:G(\bar x)<G(y)<G(\bar x)+\eta\}$, the following inequality holds:
			$$\varphi^\prime(G(x)-G(\bar x))\cdot\dist(0,\partial G(x))\ge1,$$
			where $\dist(u,S):=\inf\{\snorm{u-v}:v\in S\}$.
			\item If $G$ satisfies the K\L~property at each point of $\dom(\partial G)$, then $G$ is called a K\L~function.
		\end{enumerate}
	\end{definition}
	
	\begin{remark}\label{rem:KL}
		It is hard to tell whether the K\L~property holds for a function (at some point) by the very \Cref{def:KL property}. Nevertheless, extensive works have revealed some special cases. For example, in \cite[section 2.2]{Xu2013A2}, the authors mentioned, among others, real analytic functions, locally strongly convex functions and semi-algebraic functions. They also noted that the finite sum of real analytic and semi-algebraic functions also enjoys this property. 
	\end{remark}
	
	\par In \cite{bolte2014proximal2}, Bolte et al. proved the following uniformized version of the K\L~property, which will be used for our convergence analysis. 
	
	\begin{lemma}[Uniformized K\L~property]\label{lem:KL}
		Let $\Omega\subseteq V$ be a compact set and $G: V\to(-\infty,\infty]$ be a proper and lower semi-continuous function, where $V$ is an Euclidean space with $\inner{\cdot,\cdot}_V$ as its inner product and $\snorm{\cdot}_V=\sqrt{\langle\cdot,\cdot\rangle}_V$ as its norm. Assume that $G$ is constant on $\Omega$ and satisfies the K\L~property at each point of $\Omega$. Then, there exist $\eps>0,~\eta>0$ and $\varphi\in\Phi_\eta$ such that for all $\bar x\in\Omega$ and all $x$ in $\lrbrace{x\in V:\dist(x,\Omega)<\eps}\bigcap\{y:G(\bar x)<G(y)<G(\bar x)+\eta\}$, we have
		$$\varphi^\prime(G(x)-G(\bar x))\cdot\dist(0,\partial G(x))\ge1.$$
	\end{lemma}
	
	\section{Auxiliary Lemmas}\label{appsec:lemmas}
	In this section, we elaborate some ingredients that are useful in the later \cref{appsec:detail proof}, including
	\begin{itemize}
		\item approximate sufficient reduction of energy value at $\bar Z^{(k)}$ (\Cref{lem:approx suff reduce});
		\item approximate subgradient norm lower bound for ``iterate gap'' $\Vert \bar Z^{(k+1)}-Z^{(k)}\Vert_\F$ (\Cref{lem:approx subgrad lb}); and
		\item properties of accumulation point set of $\{\bar Z^{(k)}\}$ (\Cref{lem:properties of acc pt set}).
	\end{itemize}
	We shall mention beforehand that the  sequence $\{\bar Z^{(k)}\}$ closes the gap caused by infeasibility of the real iterate sequence $\{Z^{(k)}\}$. Thus it is not surprising that most ingredients listed above revolve around $\{\bar Z^{(k)}\}$. In the sequel, we assume $N\ge3$ by default.
	
	\par We list below constants that will be used. Note that if $\sigma>2M(N-2)$, all of them are positive.
	\begin{center}
		$\displaystyle M=\snorm{\varrho}_\infty\snorm{e}_\infty^3\snorm{C}_2+\beta,\quad C_0=\frac{\sigma}{2}-M(N-2),\quad C_1=\frac{\sigma+3MN}{2}-2M,$\\[0.2cm]
		$C_2=(N-1)C_1,\quad C_3=(\sigma+M(N-1))\sqrt{N-1}.$\\[0.2cm]
	\end{center}
	One can verify that function $f_\beta$ is $M$-Lipschitz continuously differentiable: for any $j\in\{2,\ldots,N\}$ and $Z_1,Z_2\in\big(\R^{K\times K}\big)^{N-1}$,
	\begin{equation}
		\label{eqn:lip of partial gradient}
		\snorm{\nabla_{X_j}f_\beta(Z_1)-\nabla_{X_j}f_\beta(Z_2)}_\F\le M\snorm{Z_1-Z_2}_\F.
	\end{equation}
	In addition, we use notation ``$\Delta$'' for the difference between optimal iterate \cref{eqn:optimal solution sequence restate} and real iterate. For instance, 
	\begin{equation}
		\label{eqn:difference}
		\Delta X_j^{(k+1)}:=\bar X_j^{(k+1)}-X_j^{(k+1)},\quad Z^{(k+1)}:=\bar Z^{(k+1)}-Z^{(k+1)}.
	\end{equation}
	
	\begin{lemma}\label{lem:component decrease}
		Let $\{Z^{(k)}\}$ be the sequence generated by \textup{\texttt{PBCD}}, $\{\bar Z^{(k)}\}$ be the  sequence defined in \cref{eqn:optimal solution sequence restate}. Then for each $i\in\{2,\ldots,N\}$, $\forall~ k\ge1$,
		\begin{multline}
			\label{eqn:component sufficient decrease}
			f_\beta\lrbracket{\bar X_{<i}^{(k+1)},\bar X_i^{(k)},\bar X_{>i}^{(k)}}-f_\beta\lrbracket{\bar X_{<i}^{(k+1)},\bar X_i^{(k+1)},\bar X_{>i}^{(k)}}\\\ge C_0\snorm{\bar X_i^{(k+1)}-X_i^{(k)}}_\F^2-C_1\lrbracket{\snorm{\Delta Z^{(k+1)}}_\F^2+\snorm{\Delta Z^{(k)}}_\F^2}.
		\end{multline}
	\end{lemma}
	
	\begin{proof}
		The proof mainly leverages the bilinearity of $f_\beta$ and the optimality of $\bar X_i^{(k+1)}$. We first note that the expression in the left-hand side (LHS) of \cref{eqn:component sufficient decrease} can be splitted into five parts of summations of differences below:
		\begin{enumerate}[label=\textup{Part} \arabic*]
			\item $\sum_{j=2}^{i-1}f_\beta(X_{<j}^{(k+1)},\bar X_{[j,i)}^{(k+1)},\bar X_{\ge i}^{(k)})-f_\beta(X_{\le j}^{(k+1)},\bar X_{(j,i)}^{(k+1)},\bar X_{\ge i}^{(k)})$\label{forward 1};
			\item $\sum_{j=i+1}^Nf_\beta(X_{<i}^{(k+1)},\bar X_i^{(k)},X_{(i,j)}^{(k)},\bar X_{\ge j}^{(k)})-f_\beta(X_{<i}^{(k+1)},\bar X_i^{(k)},X_{(i,j]}^{(k)},\bar X_{>j}^{(k)})$\label{forward 2};
			\item $f_\beta(X_{<i}^{(k+1)},\bar X_i^{(k)},X_{>i}^{(k)})-f_\beta(X_{<i}^{(k+1)},\bar X_i^{(k+1)},X_{>i}^{(k)})$\label{optimality};
			\item $\sum_{j=2}^{i-1}f_\beta(\bar X_{<j}^{(k+1)},X_{[j,i)}^{(k+1)},\bar X_i^{(k+1)},X_{>i}^{(k)})-f_\beta(\bar X_{\le j}^{(k+1)},X_{(j,i)}^{(k+1)},\bar X_i^{(k+1)},X_{>i}^{(k)})$\label{backward 1};
			\item $\sum_{j=i+1}^Nf_\beta(\bar X_{\le i}^{(k+1)},\bar X_{(i,j)}^{(k)},X_{\ge j}^{(k)})-f_\beta(\bar X_{\le i}^{(k+1)},\bar X_{(i,j]}^{(k)},X_{>j}^{(k)})$\label{backward 2}.
		\end{enumerate}
		
		By the optimality of $\bar X_i^{(k+1)}$ in \cref{eqn:optimal solution sequence restate} and recalling \cref{eqn:difference}, we readily have a lower bound for \cref{optimality}:
		\begin{align}
			f_\beta(X_{<i}^{(k+1)},\bar X_i^{(k)},X_{>i}^{(k)})-&f_\beta(X_{<i}^{(k+1)},\bar X_i^{(k+1)},X_{>i}^{(k)})\nonumber\\
			&\ge\frac{\sigma}{2}\lrbracket{\snorm{\bar X_i^{(k+1)}-X_i^{(k)}}_\F^2-\snorm{\Delta X_i^{(k)}}_\F^2}.\label{eqn:optimality lb}
		\end{align}
		Since the analysis for the remaining differences is very similar, we merely demonstrate in detail for \cref{forward 1,backward 1}. Recall \cref{eqn:difference},
		\begin{align*}
			\text{\cref{forward 1}}&=\sum_{j=2}^{i-1}\inner{\nabla_{X_j}f_\beta(X_{<j}^{(k+1)},\bar X_j^{(k+1)},\bar X_{(j,i)}^{(k+1)},\bar X_{\ge i}^{(k)}),\Delta X_j^{(k+1)}},\\
			\text{\cref{backward 1}}&=\sum_{j=2}^{i-1}\inner{\nabla_{X_j}f_\beta(\bar X_{<j}^{(k+1)},\bar X_j^{(k+1)},X_{(j,i)}^{(k+1)},\bar X_i^{(k+1)},X_{>i}^{(k)}),-\Delta X_j^{(k+1)}}.
		\end{align*}
		Hence using the fact that for any $a,b,c\ge0,~ab\le(a^2+b^2)/2,~(a+b)^2\le2(a^2+b^2)$ and recalling \cref{eqn:lip of partial gradient}, we have
		\begin{align}
			&\text{\cref{forward 1}}+\text{\cref{backward 1}}\label{eqn:I1+I2}\\
			\ge&-\sum_{j=2}^{i-1}M\lrbracket{\snorm{\bar X_i^{(k+1)}-\bar X_i^{(k)}}_\F+\sum_{j\ne l<i}\snorm{\Delta X_l^{(k+1)}}_\F+\sum_{l>i}\snorm{\Delta X_l^{(k)}}_\F}\snorm{\Delta X_j^{(k+1)}}_\F\nonumber\\
			\ge&-\sum_{j=2}^{i-1}\frac{M}{2}\lrbracket{3\snorm{\Delta X_j^{(k+1)}}_\F^2+\sum_{j\ne l<i}\snorm{\Delta X_l^{(k+1)}}_\F^2+\sum_{l>i}\snorm{\Delta X_l^{(k)}}_\F^2}\nonumber\\
			&-\sum_{j=2}^{i-1}M\lrbracket{\snorm{\bar X_i^{(k+1)}-X_i^{(k)}}_\F^2+\snorm{\Delta X_i^{(k)}}_\F^2}\nonumber\\
			=&-\frac{M}{2}\big[(i-3)\sum_{j=2}^{i-1}\snorm{\Delta X_j^{(k+1)}}_\F^2+(i-2)\sum_{j=i+1}^N\snorm{\Delta X_j^{(k)}}_\F^2\big]\nonumber\\
			&-\frac{3M}{2}\sum_{j=2}^{i-1}\snorm{\Delta X_j^{(k+1)}}_\F^2-M(i-2)\lrbracket{\snorm{\bar X_i^{(k+1)}-X_i^{(k)}}_\F^2+\snorm{\Delta X_i^{(k)}}_\F^2}.\nonumber
		\end{align}
		Similar arguments yield a lower bound for \cref{forward 2}$+$\cref{backward 2}:
		\begin{align}
			&\text{\cref{forward 2}}+\text{\cref{backward 2}}\label{eqn:I3+I4}\\
			\ge&-\frac{M}{2}\big[(N-i)\sum_{j=2}^{i-1}\snorm{\Delta X_j^{(k+1)}}_\F^2+(N-1-i)\sum_{j=i+1}^N\snorm{\Delta X_j^{(k)}}_\F^2\big]\nonumber\\
			&-\frac{3M}{2}\sum_{j=i+1}^N\snorm{\Delta X_j^{(k)}}_\F^2-M(N-i)\lrbracket{\snorm{\bar X_i^{(k+1)}-X_i^{(k)}}_\F^2+\snorm{\Delta X_i^{(k)}}_\F^2}\nonumber.
		\end{align}
		Combining \cref{eqn:optimality lb}, \cref{eqn:I1+I2} and \cref{eqn:I3+I4}, we have 
		\begin{align*}
			&f_\beta\lrbracket{\bar X_{<i}^{(k+1)},\bar X_i^{(k)},\bar X_{>i}^{(k)}}-f_\beta\lrbracket{\bar X_{<i}^{(k+1)},\bar X_i^{(k+1)},\bar X_{>i}^{(k)}}\\
			\ge&\left[\sigma/2-M(N-2)\right]\snorm{\bar X_i^{(k+1)}-X_i^{(k)}}_\F^2-\left[\sigma/2+M(N-2)\right]\snorm{\Delta X_i^{(k)}}_\F^2\nonumber\\
			&-\frac{3M}{2}\sum_{j<i}\snorm{\Delta X_j^{(k+1)}}_\F^2-\frac{M}{2}\big[(i-3)\sum_{j<i}\snorm{\Delta X_j^{(k+1)}}_\F^2+(i-2)\sum_{j>i}\snorm{\Delta X_j^{(k)}}_\F^2\big]\nonumber\\
			&-\frac{3M}{2}\sum_{j>i}\snorm{\Delta X_j^{(k)}}_\F^2-\frac{M}{2}\big[(N-i)\sum_{j<i}\snorm{\Delta X_j^{(k+1)}}_\F^2+(N-1-i)\sum_{j>i}\snorm{\Delta X_j^{(k)}}_\F^2\big]\nonumber\\
			=&\left[\sigma/2-M(N-2)\right]\snorm{\bar X_i^{(k+1)}-X_i^{(k)}}_\F^2-\left[\sigma/2+M(N-2)\right]\snorm{\Delta X_i^{(k)}}_\F^2\nonumber\\
			&-\frac{MN}{2}\lrbracket{\sum_{j<i}\snorm{\Delta X_j^{(k+1)}}_\F^2+\sum_{j>i}\snorm{\Delta X_j^{(k)}}_\F^2}\nonumber\\
			\ge&\left[\sigma/2-M(N-2)\right]\snorm{\bar X_i^{(k+1)}-X_i^{(k)}}_\F^2-\left[\sigma/2+M(N-2)\right]\snorm{\Delta X_i^{(k)}}_\F^2\nonumber\\
			&-\frac{MN}{2}\lrbracket{\snorm{\Delta Z^{(k+1)}}_\F^2+\snorm{\Delta Z^{(k)}}_\F^2}\nonumber,
		\end{align*}
		which completes the proof by noticing the definition of $C_0,C_1$.
	\end{proof}
	
	The following approximate sufficient reduction is then a direct corollary of the above lemma. We omit the proof.
	
	\begin{corollary}[Approximate Sufficient Reduction]\label{lem:approx suff reduce}
		Let $\{Z^{(k)}\}$ be the sequence generated by \textup{\textup{\texttt{PBCD}}}, $\{\bar Z^{(k)}\}$ be the  sequence defined in \cref{eqn:optimal solution sequence restate}. Then $\forall~ k\ge1$, 
		\begin{equation*}
			f_\beta(\bar Z^{(k)})-f_\beta(\bar Z^{(k+1)})\ge C_0\snorm{\bar Z^{(k+1)}-Z^{(k)}}_\F^2-C_2\lrbracket{\snorm{\Delta Z^{(k)}}_\F^2+\snorm{\Delta Z^{(k+1)}}_\F^2}.
		\end{equation*}
	\end{corollary}
	
	\begin{lemma}[Approximate Subgradient Lower Bound]\label{lem:approx subgrad lb}
		Let $\{Z^{(k)}\}$ be the sequence generated by \textup{\textup{\texttt{PBCD}}}, $\{\bar Z^{(k)}\}$ be the  sequence defined in \cref{eqn:optimal solution sequence restate}. Then for any $k\ge0$, there exists $W^{(k+1)}\in\partial F(\bar Z^{(k+1)})$ such that,
		\begin{equation}
			\snorm{W^{(k+1)}}_\F\le C_3\lrbracket{\snorm{\bar Z^{(k+1)}-Z^{(k)}}_\F+\snorm{\Delta Z^{(k+1)}}_\F}.
			\label{eqn:approx subgrad lb}
		\end{equation}
	\end{lemma}
	
	\begin{proof}
		For each $i\in\{2,\ldots,N\}$, it follows from the optimality of $\bar X_i^{(k+1)}$ that there exists $A_i^{(k+1)}\in \partial\delta_{\calS}(\bar X_i^{(k+1)})$ such that
		\begin{equation}
			0=\nabla_{X_i}f_\beta\lrbracket{X_{<i}^{(k+1)},\bar X_i^{(k+1)},X_{>i}^{(k)}}+\sigma\lrbracket{\bar X_i^{(k+1)}-X^{(k)}}+A_i^{(k+1)}.
			\label{eqn:element in normal cone}
		\end{equation}
		Using \cref{eqn:element in normal cone} and the calculus of subdifferential (e.g. see \cite{rockafellar2009variational2}), we have
		\begin{eqnarray*}
			\partial F(\bar Z^{(k+1)})&&\owns\nabla f_\beta(\bar Z^{(k+1)})+\lrbracket{A_i^{(k+1)}}_{i=2}^N\\
			&&=\lrbracket{\nabla_{X_i}f_\beta(\bar Z^{(k+1)})-\nabla_{X_i}f_\beta(X_{<i}^{(k+1)},\bar X_i^{(k+1)},X_{>i}^{(k)})-\sigma(\bar X_i^{(k+1)}-X_i^{(k)})}_{i=2}^N.
		\end{eqnarray*}
		
		Let $W^{(k+1)}$ be the subgradient defined above. We have for any $i\in\{2,\ldots,N\}$ 
		$$\begin{aligned}
			&\big\Vert\nabla_{X_i}f_\beta(\bar Z^{(k+1)})-\nabla_{X_i}f_\beta(X_{<i}^{(k+1)},\bar X_i^{(k+1)},X_{>i}^{(k)})-\sigma(\bar X_i^{(k+1)}-X_i^{(k)})\big\Vert_\F\\
			\le&M\lrbracket{\sum_{j<i}\snorm{\Delta X_j^{(k+1)}}_\F+\sum_{j>i}\snorm{\bar X_j^{(k+1)}-X_j^{(k)}}_\F}+\sigma\snorm{\bar X_i^{(k+1)}-X_i^{(k)}}_\F\\
			\le&M\lrbracket{\sum_{j=2}^N\snorm{\Delta X_j^{(k+1)}}_\F+\snorm{\bar X_j^{(k+1)}-X_j^{(k)}}_\F}+\sigma\snorm{\bar X_i^{(k+1)}-X_i^{(k)}}_\F\\
			\le&M\sqrt{N-1}\lrbracket{\snorm{\Delta Z^{(k+1)}}_\F+\snorm{\bar Z^{(k+1)}-Z^{(k)}}_\F}+\sigma\snorm{\bar X_i^{(k+1)}-X_i^{(k)}}_\F.
		\end{aligned}$$
		Therefore, it follows that
		$$\begin{aligned}
			\snorm{W^{(k+1)}}_\F\le&(\sigma+M(N-1))\sqrt{N-1}\snorm{\bar Z^{(k+1)}-Z^{(k)}}_\F\\
			&\qquad+M(N-1)\sqrt{N-1}\snorm{\Delta Z^{(k+1)}}_\F,
		\end{aligned}$$
		which completes the proof by recalling the definition of $C_3$.
	\end{proof}
	
	In the following lemma, we gather several properties of the accumulation point set $\omega(Z^{(0)})$ defined as
	\begin{equation}
		\omega(Z^{(0)}):=\lrbrace{\overline Z=(\overline X_i)_{i=2}^N:~\exists\text{ a subsequence }\{\bar Z^{(k)}\}_{k\in\mathcal{K}},~\st~\bar Z^{(k)}\to\overline Z\text{ in }\mathcal{K}}.
		\label{eqn:omega}
	\end{equation}
	
	\begin{lemma}[Properties of Accumulation Point Set]\label{lem:properties of acc pt set}
		Suppose $\sigma>2M(N-2)$ and \Cref{assume:iterate restate} holds. Let $\{Z^{(k)}\}$ be the sequence generated by \textup{\textup{\texttt{PBCD}}}, $\{\bar Z^{(k)}\}$ be the  sequence defined in \cref{eqn:optimal solution sequence restate}. Then we have that
		\begin{enumerate}[label=\textup{(}\arabic*\textup{)}]
			\item $\omega(Z^{(0)})\subseteq\calS^{N-1}$ is nonempty compact and $\dist(\bar Z^{(k)},\omega(Z^{(0)}))\to0$;\label{lem:acc:compact}
			\item $F$ is finite and constant on $\omega(Z^{(0)})$;\label{lem:acc:finite and constant} 
			\item all the elements in $\omega(Z^{(0)})$ are KKT points of \cref{eqn:penalty problem restate}.\label{lem:acc:KKT}
		\end{enumerate}
	\end{lemma}
	
	\begin{proof}
		(1) Since $\{\bar Z^{(k)}\}\subseteq\calS^{N-1}$ and $\calS$ is bounded, there exist a subsequence $\{\bar Z^{(k)}\}_{k\in\mathcal{K}}$ and $\bbar Z$ such that $\bar Z^{(k)}\to \bbar Z$ in $\mathcal{K}$, which gives $\omega(Z^{(0)})\ne\emptyset$. Also note that $\calS^{N-1}$ is closed, thus $\bbar Z\in\calS^{N-1}$ and hence $\omega(Z^{(0)})\subseteq\calS^{N-1}$. To prove the compactness, we only need to verify the closedness. Let $\omega(Z^{(0)})\owns\bbar Z_j\to\widetilde Z$. For any $j$, there exists a subsequence $\mathcal{K}_j\subseteq\mathbb{N}$ so that $\bar Z^{(k)}\to\bbar Z_j$ in $\mathcal{K}_j$. Hence for each $j$, we can pick an integer $k_j\in\mathcal{K}_j$ such that $\{k_j\}_j$ is a sequence increasing to $\infty$, and $\snorm{\bar Z^{(k_j)}-\bbar Z_j}_\F\le\frac{1}{j}$. Therefore, the sequence $\{\bar Z^{(k_j)}\}_j$ satisfies
		$$\snorm{\bar Z^{(k_j)}-\widetilde Z}_\F\le\snorm{\bar Z^{(k_j)}-\bbar Z_j}_\F+\snorm{\bbar Z_j-\widetilde Z}_\F\le\frac{1}{j}+\snorm{\bbar Z_j-\widetilde Z}_\F.$$
		Based on $\bbar Z_j\to\widetilde Z$ and the definition \cref{eqn:omega} of $\omega(Z^{(0)})$, we obtain $\widetilde Z\in\omega(Z^{(0)})$, and hence the closedness of $\omega(Z^{(0)})$. 
		\par The latter part of this item follows directly from the definition \cref{eqn:omega} of the nonempty $\omega(Z^{(0)})$.
		
		\vskip 0.1cm
		
		\par (2) By (1), $F$ must be finite on $\omega(Z^{(0)})\subseteq\calS^{N-1}$. Now for any $\overline Z\in\omega(Z^{(0)})$, there exists a subsequence $\{\bar Z^{(k)}\}_{k\in\mathcal{K}}$ converging to $\overline Z$. The monotonicity of $\{F(\bar Z^{(k)})\}$ (\Cref{assume:iterate restate} \cref{assume:decreasing func 2}) and the continuity of $F$ in $\calS^{N-1}$ imply $F(\bbar Z) = \bbar F$. By the arbitrariness of $\overline Z$, $F\equiv \bbar F$ on $\omega(Z^{(0)})$.
		
		\vskip 0.1cm
		\par (3) From \Cref{lem:approx suff reduce}, we get for $\forall~ k\ge1$,
		$$C_0\snorm{\bar Z^{(k+1)}-Z^{(k)}}_\F^2\le f_\beta(\bar Z^{(k)})-f_\beta(\bar Z^{(k+1)})+C_2\lrbracket{\snorm{\Delta Z^{(k)}}_\F^2+\snorm{\Delta Z^{(k+1)}}_\F^2}.$$
		Summing the above inequality over $k$ from $r$ to $s$ ($s\ge r>1$), we have
		\begin{eqnarray}
			C_0\sum_{k=r}^s\snorm{\bar Z^{(k+1)}-Z^{(k)}}_\F^2\le&&\sum_{k=r}^s\lrbracket{f_\beta(\bar Z^{(k)})-f_\beta(\bar Z^{(k+1)})}+2C_2\sum_{k=r-1}^{s}\snorm{\Delta Z^{(k)}}_\F^2\nonumber\\
			=&&f_\beta(\bar Z^{(r)})-f_\beta(\bar Z^{(s+1)})+2C_2\sum_{k=r-1}^{s}\snorm{\Delta Z^{(k)}}_\F^2.
			\label{eqn:iter gap to 0}
		\end{eqnarray}
		\Cref{assume:iterate restate} \cref{assume:decreasing func 2} and $\sigma>2M(N-2)$ (hence $C_0>0$) imply 
		$$\sum_{k=r}^s\snorm{\bar Z^{(k+1)}-Z^{(k)}}_\F^2\le C_0^{-1}\lrbracket{f_\beta(\bar Z^{(r)})-\bbar F+2C_2\sum_{k=1}^\infty\Vert \Delta Z^{(k)}\Vert_\F^2}<\infty,~\forall~ s\ge r>1.$$
		Therefore $\Vert \bar Z^{(k+1)}-Z^{(k)}\Vert_\F\to0$. By \Cref{lem:approx subgrad lb} and \Cref{assume:iterate restate} \cref{assume:local solver 2}, we further derive $W^{(k+1)}\to0$.
		
		\par Now pick $\overline Z\in\omega(Z^{(0)})$ and the associated converging subsequence $\{\bar Z^{(k+1)}\}_{k\in\mathcal{K}}$. Since $\{(\bar Z^{(k+1)},W^{(k+1)})\}_{k\in\mathcal{K}}$ converges to $\lrbracket{\bbar Z,0}$, $(\bar Z^{(k+1)},W^{(k+1)})\in\mathrm{graph}(\partial F)$ and $F(\bar Z^{(k+1)})\to \bbar F=F(\bbar Z)$ (by (2)), from \Cref{rem:subdiff} \cref{rem:closedness of graph}, we deduce $(\bbar Z,0)\in\mathrm{graph}(\partial F)$, or equivalently, $0\in\partial F(\bbar Z)$. In view of \Cref{rem:subdiff} \cref{rem:critical and KKT}, $\bbar Z$ is a KKT point of \cref{eqn:penalty problem restate}. We complete the proof by the arbitrariness of $\bbar Z$. 
	\end{proof}
	
	\section{Formal Statement and Proof of Theorem 3.3}\label{appsec:detail proof}
	%\section{Formal Statement and Proof of Theorem \ref{thm:rough convergence property}}\label{appsec:detail proof}
	
	Thus far, we have established three auxiliary ingredients and are ready to prove Theorem 3.3 using the uniformized K\L~property (\Cref{lem:KL}). Beforehand, we shall mention that $f_\beta$ is a real analytic function, and $\delta_\calS$ is semi-algebraic since $\calS$ is apparently a semi-algebraic set. According to \Cref{rem:KL}, $F$ is a K\L~function, which is the premise for applying \Cref{lem:KL}. 
	
	\par For a global minimizer $\bbar Z$, since $F$ is a K\L~function, we have the related scalar $\bar\eta>0$, neighborhood $\bbar{\mathcal{U}}$ of $\bbar Z$ and function $\bar\varphi\in\Phi_{\bar\eta}$ as in \Cref{def:KL property}, which will be used in the sequel. Additionally, we define the following constants. 
	\vskip 0.1cm
	\begin{center}
		$\displaystyle\tilde C_0=\sqrt{C_0}-\frac{1}{2p},\quad\tilde C_1=\frac{pC_3}{\tilde C_0},$\\[0.2cm]
		$\displaystyle\tilde C_2=\lrbracket{\frac{1}{2p\tilde C_0}+1}\frac{1}{\sqrt{C_0}},\quad\tilde C_3=\frac{1}{\tilde C_0}\lrbracket{\frac{1}{2p}+2\sqrt{C_2}}+\sqrt{C_2}\tilde C_2+1.$
	\end{center}
	\vskip 0.2cm
	Here constant $p$ satisfies 
	\begin{equation}p>\frac{1}{2\sqrt{C_0}}.\label{eqn:p}\end{equation}
	As a result, all the constants listed above are positive.
	
	\par In what follows, we give the formal statement of Theorem 3.3 and its complete proof.
	
	\begin{theorem}[Convergence Property of Algorithm 2.3 -- Inexact Version]
		Suppose $\sigma>2M(N-2)$ and \Cref{assume:iterate restate} holds. Let $\{Z^{(k)}\}$ be the sequence generated by Algorithm 2.3, $\{\bar Z^{(k)}\}$ be the  sequence defined in \cref{eqn:optimal solution sequence restate}. Assume also that $F(\bar Z^{(k)})>\bbar F,~\forall~ k\ge1$. Then 
		\begin{enumerate}[label=\textup{(}\arabic*\textup{)}]
			\item the sequence $\{Z^{(k)}\}$ converges to a KKT point of \cref{eqn:penalty problem restate};
			\item if further $Z^{(0)}\in\{Z:F(\bbar Z)<F(Z)<F(\bbar Z)+\bar\eta\}$ and
			\begin{align}
				&\tilde C_1\bar\varphi\left(F(Z^{(0)})-F(\bbar Z)\right)+\tilde C_2\sqrt{F(Z^{(0)})-F(\bbar Z)}\nonumber\\
				&\qquad\qquad\qquad+\tilde C_3\sum_{k=1}^\infty\snorm{\Delta Z^{(k)}}_\F+\snorm{Z^{(0)}-\bbar Z}_\F\le r,
				\label{eqn:good initial pt}
			\end{align}
			where $\bbar Z$ is a global minimizer of \cref{eqn:penalty problem restate} and $r$ is a positive scalar such that the closed ball $B_r(\bbar Z):=\{Z:\snorm{Z-\bbar Z}_\F\le r\}\subseteq\bbar{\mathcal{U}}$. Then $\{Z^{(k)}\}$ converges to a global minimizer of \cref{eqn:penalty problem restate}. 
		\end{enumerate}
		\label{thm:convergence inexact}
	\end{theorem}
	
	\begin{proof}
		(1) Let $Z_*\in\omega(Z^{(0)})$ and $\varphi\in\Phi_\eta,~\eps,\eta>0$ be defined in \Cref{lem:KL}, where $G:=F,~\Omega:=\omega(Z^{(0)})$. Since $F(\bar Z^{(k)})>\bbar F=F(Z_*)$ and $\{F(\bar Z^{(k)})\}$ converges monotonically to $\bbar F$ (\Cref{assume:iterate restate} \cref{assume:decreasing func 2}), there exists an integer $k_0$, such that $F(Z_*)<F(\bar Z^{(k)})<F(Z_*)+\eta,~\forall~ k\ge k_0$. Then $\varphi^\prime(F(\bar Z^{(k)})-F(Z_*))$ is well-defined for any $k\ge k_0$. On the other, by \Cref{lem:properties of acc pt set} \cref{lem:acc:compact}, there exists an integer $k_1$, such that $\dist(\bar Z^{(k)},\omega(Z^{(0)}))<\eps,~\forall~ k\ge k_1$. Let $s=\max(k_0,k_1)$. By \Cref{lem:KL}, for $\forall~ k\ge s$,
		$$\varphi^\prime\lrbracket{F(\bar Z^{(k)})-F(Z_*)}\cdot\dist\lrbracket{0,\partial F(\bar Z^{(k)})}\ge1,$$
		which, combined with \Cref{lem:approx subgrad lb}, further gives
		\begin{align}
			\varphi^\prime\lrbracket{F(\bar Z^{(k)})-F(Z_*)}&\ge\frac{1}{\dist(0,\partial F(\bar Z^{(k)}))}\ge\frac{1}{\Vert W^{(k)}\Vert_\F}\nonumber\\
			&\ge\frac{1}{C_3(\Vert \bar Z^{(k)}-Z^{(k-1)}\Vert_\F+\Vert \Delta Z^{(k)}\Vert_\F)}.
			\label{eqn:KL}
		\end{align}
		
		Let 
		$$\D^{(k,k+1)}:=\varphi\lrbracket{F(\bar Z^{(k)})-F(Z_*)}-\varphi\lrbracket{F(\bar Z^{(k+1)})-F(Z_*)}.$$ 
		By the concavity of $\varphi$, \Cref{assume:iterate restate} \cref{assume:decreasing func 2}, \Cref{lem:approx suff reduce} and \cref{eqn:KL}, we obtain
		\begin{align*}
			\D^{(k,k+1)}&\ge\varphi^\prime\lrbracket{F(\bar Z^{(k)})-F(Z_*)}\cdot\lrbracket{F(\bar Z^{(k)})-F(\bar Z^{(k+1)})}\\
			&\ge\frac{F(\bar Z^{(k)})-F(\bar Z^{(k+1)})}{C3(\Vert \bar Z^{(k)}-Z^{(k-1)}\Vert_\F+\Vert \Delta Z^{(k)}\Vert_\F)}\\
			&\ge\frac{C_0\Vert \bar Z^{(k+1)}-Z^{(k)}\Vert_\F^2-C_2(\Vert \Delta Z^{(k)}\Vert_\F^2+\Vert \Delta Z^{(k+1)}\Vert_\F^2)}{C_3(\Vert \bar Z^{(k)}-Z^{(k-1)}\Vert_\F+\Vert \Delta Z^{(k)}\Vert_\F)}.
		\end{align*}
		An upper bound for $\Vert \bar Z^{(k+1)}-Z^{(k)}\Vert_\F$ follows from some algebraic manipulations:
		\begin{align}
			\sqrt{C_0}\snorm{\bar Z^{(k+1)}-Z^{(k)}}_\F\le&\sqrt{C_3\D^{(k,k+1)}\Vert \bar Z^{(k)}-Z^{(k-1)}\Vert_\F}+\sqrt{C_3\D^{(k,k+1)}\Vert \Delta Z^{(k)}\Vert_\F}\nonumber\\
			&+\sqrt{C_2}\snorm{\Delta Z^{(k)}}_\F+\sqrt{C_2}\snorm{\Delta Z^{(k+1)}}_\F\nonumber\\
			\le&pC_3\D^{(k,k+1)}+\frac{1}{2p}\lrbracket{\Vert \bar Z^{(k)}-Z^{(k-1)}\Vert_\F+\Vert \Delta Z^{(k)}\Vert_\F}\label{eqn:iter gap ub}\\
			&+\sqrt{C_2}\snorm{\Delta Z^{(k)}}_\F+\sqrt{C_2}\snorm{\Delta Z^{(k+1)}}_\F\nonumber,
		\end{align}
		where the first inequality uses $\sqrt{a+b}\le\sqrt{a}+\sqrt{b},~\forall~ a,b\ge0$, while the second one uses $\sqrt{ab}\le pa/2+b/(2p),~\forall~ a,b\ge0,~p>0$. Here $p$ is chosen as in \cref{eqn:p}. Summing \cref{eqn:iter gap ub} over $k$ from $s$ to some $t\ge s$ gives
		\begin{align}
			&\sqrt{C_0}\sum_{k=s}^t\snorm{\bar Z^{(k+1)}-Z^{(k)}}_\F\nonumber\\
			\le&pC_3\sum_{k=s}^t\D^{(k,k+1)}+\frac{1}{2p}\lrbracket{\sum_{k=s}^t\Vert \bar Z^{(k+1)}-Z^{(k)}\Vert_\F+\Vert \bar Z^{(s)}-Z^{(s-1)}\Vert_\F}\nonumber\\
			&+\lrbracket{\frac{1}{2p}+\sqrt{C_2}}\sum_{k=s}^{t}\snorm{\Delta Z^{(k)}}_\F+\sqrt{C_2}\sum_{k=s+1}^{t+1}\snorm{\Delta Z^{(k)}}_\F\nonumber\\
			\le&pC_3\D^{(s,t+1)}+\frac{1}{2p}\lrbracket{\sum_{k=s}^t\Vert \bar Z^{(k+1)}-Z^{(k)}\Vert_\F+\Vert \bar Z^{(s)}-Z^{(s-1)}\Vert_\F}\nonumber\\
			&+\lrbracket{\frac{1}{2p}+2\sqrt{C_2}}\sum_{k=s}^{t+1}\snorm{\Delta Z^{(k)}}_\F.
		\end{align}
		Since
		$$\D^{(s,t+1)}=\varphi\lrbracket{F(\bar Z^{(s)})-F(Z_*)}-\varphi\lrbracket{F(\bar Z^{(t+1)})-F(Z_*)}\le\varphi\lrbracket{F(\bar Z^{(s)})-F(Z_*)},$$
		we further derive that
		\begin{align}
			\lrbracket{\sqrt{C_0}-\frac{1}{2p}}\sum_{k=s}^t\snorm{\bar Z^{(k+1)}-Z^{(k)}}_\F\le&pC_3\varphi\big(F(\bar Z^{(s)})-F(Z_*)\big)+\frac{1}{2p}\snorm{\bar Z^{(s)}-Z^{(s-1)}}_\F\nonumber\\
			&+\lrbracket{\frac{1}{2p}+2\sqrt{C_2}}\sum_{k=s}^{t+1}\snorm{\Delta Z^{(k)}}_\F.
			\label{eqn:iter gap ub 2}
		\end{align}
		Passing $t\to\infty$ in \cref{eqn:iter gap ub 2}, and recalling \cref{eqn:p} and the summability of $\{\Vert \Delta Z^{(k)}\Vert_\F\}$ (\Cref{assume:iterate restate} \cref{assume:local solver 2}), we conclude that $\{\Vert \bar Z^{(k+1)}-Z^{(k)}\Vert_\F\}$ is summable. Again using the summability of $\{\Vert \Delta Z^{(k)}\Vert_\F\}$, we have the convergence of $\{Z^{(k)}\}$: for any $t\ge s$, 
		$$\sum_{k=s}^t\snorm{Z^{(k+1)}-Z^{(k)}}\le\sum_{k=s}^t\snorm{\Delta Z^{(k+1)}}_\F+\sum_{k=s}^t\snorm{\bar Z^{(k+1)}-Z^{(k)}}_\F.$$
		By \Cref{lem:properties of acc pt set} \cref{lem:acc:KKT}, $Z_*\in\omega(Z^{(0)})$ is a KKT point of \cref{eqn:penalty problem restate}. Since $\Vert\Delta Z^{(k)}\Vert_\F\to0$, $\{\bar Z^{(k)}\}$ and $\{Z^{(k)}\}$ converge to the same point. Therefore, $\omega(Z^{(0)})=\{Z_*\}$ and $Z^{(k)}\to Z_*$. 
		
		\vskip 0.2cm 
		
		(2) We first assume the following claim holds: 
		\begin{claim}\label{clm:basin}
			Under the assumptions made in (2), we have 
			$$\bar Z^{(k)}\in B_r(\bbar Z)\bigcap\{Z:F(\bbar Z)<F(Z)<F(\bbar Z)+\bar\eta\},\quad\forall~ k\ge0.$$
		\end{claim}
		By (1), we derive that $\{\bar Z^{(k)}\}$ (and hence $\{Z^{(k)}\}$) converges to some $Z_*\in B_r(\bbar Z)$ and $0\in\partial F(Z_*)$. If $F(Z_*)>F(\bbar Z)$, then the K\L~property at $\bbar Z$ indicates
		$$\bar\varphi^\prime\lrbracket{F(Z_*)-F(\bbar Z)}\cdot\dist\lrbracket{0,\partial F(Z_*)}\ge1,$$
		which contradicts with $0\in\partial F(Z_*)$. Therefore, we must have $F(Z_*)=F(\bbar Z)$. 
		
		\par Now We turn to prove \Cref{clm:basin} by induction. Denote $\bar Z^{(0)}=Z^{(0)}$. The \Cref{clm:basin} holds trivially for $k=0$. Next let $k=1$. Since $Z^{(0)}$ is feasible, by \Cref{lem:approx suff reduce}, we have 
		\begin{equation*}
			F\lrbracket{Z^{(0)}}-F\lrbracket{\bar Z^{(1)}}=F\lrbracket{\bar Z^{(0)}}-F\lrbracket{\bar Z^{(1)}}\ge C_0\snorm{\bar Z^{(1)}-Z^{(0)}}_\F^2-C_2\snorm{\Delta Z^{(1)}}_\F^2,
		\end{equation*}
		which, together with $\sigma>2M(N-2)$, yield
		\begin{align}
			\snorm{\bar Z^{(1)}-Z^{(0)}}_\F&\le\frac{1}{\sqrt{C_0}}\sqrt{(F(Z^{(0)})-F(\bar Z^{(1)}))+C_2\Vert \Delta Z^{(1)}\Vert_\F^2}\nonumber\\
			&\le\frac{1}{\sqrt{C_0}}\lrbracket{\sqrt{F(Z^{(0)})-F(\bar Z^{(1)})}+\sqrt{C_2}\Vert \Delta Z^{(1)}\Vert_\F}\label{eqn:induction k=1}\\
			&\le\frac{1}{\sqrt{C_0}}\lrbracket{\sqrt{F\lrbracket{Z^{(0)}}-F(\bbar Z)}+\sqrt{C_2}\Vert \Delta Z^{(1)}\Vert_\F}.\nonumber
		\end{align}
		Combining \cref{eqn:induction k=1}, triangle inequality of Frobenius norm and \cref{eqn:good initial pt}, we deduce
		\begin{align*}
			\snorm{\bar Z^{(1)}-\bbar Z}_\F&\le\snorm{\bar Z^{(1)}-Z^{(0)}}_\F+\snorm{Z^{(0)}-\bbar Z}_\F\\
			&\le\frac{1}{\sqrt{C_0}}\sqrt{F\lrbracket{Z^{(0)}}-F\lrbracket{\bbar Z}}+\sqrt{\frac{C_2}{C_0}}\snorm{\Delta Z^{(1)}}_\F+\snorm{Z^{(0)}-\bbar Z}_\F\le r.
		\end{align*}
		Hence $\bar Z^{(1)}\in B_r(\bbar Z)$, and \Cref{clm:basin} holds for $k=1$ by noting $F(\bbar Z)<F(\bar Z^{(1)})\le F(\bar Z^{(0)})=F(Z^{(0)})<F(\bbar Z)+\bar\eta$. 
		
		\par Now suppose \Cref{clm:basin} holds for some $t\ge1$. From \cref{eqn:iter gap ub 2} in the proof of (1), we get 
		\begin{align}
			\tilde C_0\sum_{k=1}^t\snorm{\bar Z^{(k+1)}-Z^{(k)}}_\F\le&pC_3\bar\varphi\lrbracket{F(\bar Z^{(1)})-F(\bbar Z)}+\frac{1}{2p}\snorm{\bar Z^{(1)}-Z^{(0)}}_\F\nonumber\\
			&+\lrbracket{\frac{1}{2p}+2\sqrt{C_2}}\sum_{k=1}^{t+1}\snorm{\Delta Z^{(k)}}_\F\nonumber\\
			\le&pC_3\bar\varphi\lrbracket{F(Z^{(0)})-F(\bbar Z)}+\frac{1}{2p}\snorm{\bar Z^{(1)}-Z^{(0)}}_\F\label{eqn:iter gap ub 2 s=1}\\
			&+\lrbracket{\frac{1}{2p}+2\sqrt{C_2}}\sum_{k=1}^{t+1}\snorm{\Delta Z^{(k)}}_\F,\nonumber
		\end{align}
		where the second inequality comes from the monotonicity of $\{F(\bar Z^{(k)})\}$ (\Cref{assume:iterate restate} \cref{assume:decreasing func 2}) and $\bar\varphi^\prime>0$. Combining \cref{eqn:iter gap ub 2 s=1} and, again, \cref{eqn:good initial pt}, \cref{eqn:induction k=1} and triangle inequality, we derive that 
		\begin{align*}
			\snorm{\bar Z^{(t+1)}-\bbar Z}_\F\le&\sum_{k=1}^t(\snorm{\bar Z^{(k+1)}-Z^{(k)}}_\F+\snorm{\Delta Z^{(k)}}_\F)+\snorm{\bar Z^{(1)}-Z^{(0)}}_\F+\snorm{Z^{(0)}-\bbar Z}_\F\\
			\le&\frac{pC_3}{\tilde C_0}\bar\varphi(F(Z^{(0)})-F(\bbar Z))+\frac{1}{\tilde C_0}\lrbracket{\frac{1}{2p}+2\sqrt{C_2}}\sum_{k=1}^{t+1}\snorm{\Delta Z^{(k)}}_\F\\
			&+\frac{\snorm{\bar Z^{(1)}-Z^{(0)}}_\F}{2p\tilde C_0}+\sum_{k=1}^t\snorm{\Delta Z^{(k)}}_\F+\snorm{\bar Z^{(1)}-Z^{(0)}}_\F+\snorm{Z^{(0)}-\bbar Z}_\F\\
			\le&\frac{pC_3}{\tilde C_0}\bar\varphi\lrbracket{F(Z^{(0)})-F(\bbar Z)}+\lrbracket{\frac{1}{2p\tilde C_0}+1}\snorm{\bar Z^{(1)}-Z^{(0)}}_\F\\
			&+\big[\frac{1}{\tilde C_0}\lrbracket{\frac{1}{2p}+2\sqrt{C_2}}+1\big]\sum_{k=1}^{t+1}\snorm{\Delta Z^{(k)}}_\F+\snorm{Z^{(0)}-\bbar Z}_\F\\
			\le&\tilde C_1\bar\varphi\lrbracket{F(Z^{(0)})-F(\bbar Z)}+\tilde C_2\sqrt{F(Z^{(0)})-F(\bbar Z)}\\
			&+\tilde C_3\sum_{k=1}^\infty\snorm{\Delta Z^{(k)}}_\F+ \snorm{Z^{(0)}-\bbar Z}_\F\le r.
		\end{align*}
		Hence $\bar Z^{(t+1)}\in B_r(\bbar Z)$, and \Cref{clm:basin} holds for $t+1$ by noting $F(\bbar Z)<F(\bar Z^{(t+1)})\le F(\bar Z^{(0)})=F(Z^{(0)})<F(\bbar Z)+\bar\eta$. By induction, \Cref{clm:basin} holds for any $k\ge0$. We complete the proof.
	\end{proof}
	
	\begin{remark}
		\begin{enumerate}[label=\textup{(}\arabic*\textup{)}]
			\item The results in this section recover those established in \cite{bolte2014proximal2} if $Z^{(k)}=\bar Z^{(k)},~\forall~ k$.
			\item The results in this section still hold if $f$ is Lipschitz continuously differentiable and \textup{\texttt{PBCD}} performs proximal linearized minimization for each block. 
		\end{enumerate}
	\end{remark}
	
	\section{Numerical Comparison among Local Solvers}\label{subsec:comparison NLP}
	
	Since few methods are known for MPGCC, we consider (1.8) with $N=3$ (hence an MPCC), and compare our local algorithm \textup{\texttt{PBCD}} with some modified nonlinear programming solvers, including \texttt{Knitro} \cite{byrd2006knitro2,waltz2006interior2} and \texttt{filtermpec} \cite{fletcher2004solving2,fletcher2006local2}. Version 12.4.0 of \texttt{Knitro} is available in the downloadable AMPL system \cite{fourer2003ampl2}, and \texttt{filtermpec} has an interface to the online NEOS Server \cite{czyzyk1998neos2,dolan2001neos2,gropp1997neos2}.
	
	\par We investigate the performance of \texttt{Knitro}, \texttt{filtermpec} and \textup{\texttt{PBCD}} on three 1D $N=3$ systems in subsection 4.2, starting from the initial points provided by the GR subroutine (see Table 4.2 (a)). Additionally in virtue of the high-quality initial points, we set \texttt{alg} = 4 for \texttt{Knitro} so as to invoke sequential quadratic programming algorithm\footnote{For reasons, please refer to \url{https://www.artelys.com/docs/knitro//2_userGuide/algorithms.html\#sec-algorithms}.}. We use up to 40 cores for computation involving \texttt{Knitro} and \textup{\texttt{PBCD}} and use AMPL input for \texttt{Knitro} and \texttt{filtermpec}. We monitor the repulsive energy (denoted by E), the average error (denoted by err) and the wall clock time in seconds (denoted by CPU) of all the three solvers. Note the CPU time reported by \texttt{filtermpec} is based on the computational resources on the NEOS Server. We need to point out that the formulation fed to \texttt{Knitro} and \texttt{filtermpec} is the original MPCC (1.8) rather than the $\ell_1$ penalty form \cref{eqn:penalty problem restate}.
	
	\par We terminate these solvers according to optimality violation. For \textup{\texttt{PBCD}}, the sub-solver \texttt{SSNCG} stops if primal infeasibility (defined in (4.2)) is smaller than $\eps_{\mathrm{inner}}=10^{-9}$, and the outer loop terminates when KKT violation\footnote{The expression of KKT violation can be derived by comparing the iterative scheme of \textup{\texttt{PBCD}} \cref{eqn:PBCD scheme} and KKT condition in \Cref{lem:KKT}. Note we actually omit $\mathrm{Feas}(X_i^{(k)})$ in \cref{eqn:KKT violation} due to its tiny magnitude.} 
	\begin{equation}\label{eqn:KKT violation}\sum_{i=2}^N\big\Vert\sum_{j>i}[\Xi\Lambda(X_j^{(k-1)}-X_j^{(k)})\Xi C\Xi+\beta(X_j^{(k-1)}-X_j^{(k)})]+\sigma(X_i^{(k)}-X_i^{(k-1)})\big\Vert_\F\end{equation}
	is smaller than a prescribed value $\eps_{\mathrm{outer}}$, chosen as in (4.3). In \texttt{Knitro}, we specify stopping criteria \texttt{feastol\_abs} $=$ \texttt{opttol\_abs} $=\eps_{\mathrm{outer}}$\footnote{Actually, \texttt{Knitro} stops as well in default setting if stepsize is smaller than $10^{-15}$ for 3 consecutive iterations. We keep this for less running time of \texttt{Knitro}.}; in \texttt{filtermpec}, we set \texttt{eps} $=\eps_{\mathrm{outer}}$. Other parameters of these two solvers are left unchanged as default.
	
	\par The results on three 1D $N=3$ examples are gathered in \Cref{tab:compare NLP}, where the best energies, average errors and CPUs are marked out in bold. The average errors in the brackets are the those of the corresponding initial points (see also the ``err\_s'' column in Table 4.2 (a)). 
	
	\begin{table}[htbp]
		\footnotesize
		\centering
		\caption{Comparison among local solvers}
		\label{tab:compare NLP}
		\resizebox{\textwidth}{32mm}{
			\begin{tabular}{l||rrr||rrr||rrr||rrr}
				\toprule
				\multicolumn{13}{c}{Example 1}\\\midrule[0.02em]
				\multirow{2}{*}{Alg.} & \multicolumn{3}{c||}{$K=24$ (err = $0.049$)} & \multicolumn{3}{c||}{$K=48$ (err = $0.022$)} & \multicolumn{3}{c||}{$K=96$ (err = $0.014$)} & \multicolumn{3}{c}{$K=192$ (err = $0.007$)} \\\cmidrule{2-13}
				& E & err & CPU & E & err & CPU & E & err & CPU & E & err & CPU \\\midrule[0.02em]
				\texttt{Knitro} & \textbf{18.911} & \textbf{0.013} & 4.11 & \textbf{19.004} & \textbf{0.009} & 9.33 & \textbf{19.019} & \textbf{0.004} & 393.41 & \textbf{19.021} & \textbf{0.003} & 7487.25 \\
				\texttt{filtermpec} & 19.035 & 0.071 & 12.65 & \textbf{19.004} & \textbf{0.009} & 2.74 & \textbf{19.019} & \textbf{0.004} & 29.90 & \textbf{19.021} & \textbf{0.003} & 531.72 \\
				\textup{\texttt{PBCD}} & \textbf{18.911} & \textbf{0.013} & \textbf{0.55} & \textbf{19.004} & \textbf{0.009} & \textbf{1.67} & \textbf{19.019} & \textbf{0.004} & \textbf{8.83} & \textbf{19.021} & \textbf{0.003} & \textbf{110.58} \\\midrule[0.07em]
				\multicolumn{13}{c}{Example 2}\\\midrule[0.02em]

				\multirow{2}{*}{Alg.} & \multicolumn{3}{c||}{$K=24$ (err = $0.041$)} & \multicolumn{3}{c||}{$K=48$ (err = $0.026$)} & \multicolumn{3}{c||}{$K=96$ (err = $0.017$)} & \multicolumn{3}{c}{$K=192$ (err = $0.010$)} \\\cmidrule{2-13}
				& E & err & CPU & E & err & CPU & E & err & CPU & E & err & CPU \\\midrule[0.02em]
				\texttt{Knitro} & \textbf{12.372} & 0.013 & 5.30 & \textbf{12.367} & \textbf{0.012} & 1.36 & \textbf{12.360} & \textbf{0.009} & 18.37 & \textbf{12.358} & \textbf{0.003} & 11199.66 \\
				\texttt{filtermpec} & 12.465 & 0.069 & 7.06 & 12.368 & 0.015 & 2.40 & \textbf{12.360} & \textbf{0.009} & 30.54 & \textbf{12.358} & 0.004 & 523.86 \\
				\textup{\texttt{PBCD}} & \textbf{12.372} & \textbf{0.011} & \textbf{0.41} & \textbf{12.367} & \textbf{0.012} & \textbf{0.59} & \textbf{12.360} & \textbf{0.009} & \textbf{3.16} & \textbf{12.358} & \textbf{0.003} & \textbf{25.82} \\\midrule[0.07em]

				\multicolumn{13}{c}{Example 3}\\\midrule[0.02em]
				\multirow{2}{*}{Alg.} & \multicolumn{3}{c||}{$K=24$ (err = $0.053$)} & \multicolumn{3}{c||}{$K=48$ (err = $0.027$)} & \multicolumn{3}{c||}{$K=96$ (err = $0.026$)} & \multicolumn{3}{c}{$K=192$ (err = $0.013$)} \\\cmidrule{2-13}
				& E & err & CPU & E & err & CPU & E & err & CPU & E & err & CPU \\\midrule[0.02em]
				\texttt{Knitro} & \textbf{6.318} & \textbf{0.018} & 5.53 & \textbf{6.389} & \textbf{0.013} & 17.53 & \textbf{6.400} & \textbf{0.011} & 178.23 & \textbf{6.403} & \textbf{0.001} & 19890.34 \\
				\texttt{filtermpec} & 6.383 & 0.105 & 12.00 & \textbf{6.389} & \textbf{0.013} & 2.39 & \textbf{6.400} & \textbf{0.011} & 37.94 & \textbf{6.403} & \textbf{0.001} & 504.53 \\
				\textup{\texttt{PBCD}} & \textbf{6.318} & \textbf{0.018} & \textbf{0.41} & \textbf{6.389} & \textbf{0.013} & \textbf{1.97} & \textbf{6.400} & \textbf{0.011} & \textbf{3.63} & \textbf{6.403} & \textbf{0.001} & \textbf{58.19} \\\bottomrule
		\end{tabular}}
	\end{table}
	
	\par The results in \Cref{tab:compare NLP} demonstrate that \texttt{PBCD} is comparable to both \texttt{Knitro} and \texttt{filtermpec}, the two state-of-art solvers for MPCC, in view of converged energies and solution errors, and meanwhile, underscore the incredible superiority of \textup{\texttt{PBCD}} over the other two in terms of CPU time. This is not surprising, because both of them need to handle an $O(K^4)$ sparse Hessian, whose number of non-zeros is of $O(K^3)$. Indeed though not listed, the online interface of \texttt{filtermpec} fails to handle (1.8) with size $K=384$ owing to the explosion of memory. Moreover according to our numerical experience, the proximal parameter $\sigma$ can be tuned for better performance of \texttt{PBCD}. To conclude, our proposed \texttt{PBCD} is a nice alternative in solving (1.8) when $N=3$, and emerges as a novel and efficient tool for $N>3$ cases. 
	
	\normalem
	\bibliographystyle{siamplain}
	\bibliography{sm_note}